\documentclass{llncs}
\usepackage{makeidx}  
\usepackage[pdftex]{graphicx}
\usepackage{amsmath}
\usepackage{amsfonts,amssymb}
\usepackage{clrscode3e}            
\usepackage{float}

\begin{document}
\frontmatter          
\pagestyle{headings}  
\addtocmark{Mixed-precision CG} 
\title{Impact of metric characteristics\\ of the matrix graph on the convergence\\ of the conjugate gradient method}
\titlerunning{Impact of Metric Characteristics of the Matrix Graph on the Convergence of the CG method}  
%
\author{Alexander Prolubnikov}
\authorrunning{A. Prolubnikov} 
%
\tocauthor{Alexander Prolubnikov}
\institute{Novosibirsk State University\\
\email{a.v.prolubnikov@mail.ru}
}
\maketitle              

\begin{abstract}
A model of propagation and accumulation of errors caused by round-off errors during conjugate gradient (CG) iterations is proposed, based on their interpretation as iterations of traversals of the linear system matrix graph. In the matrix graph, vertices correspond to the components of the search direction, residual, and approximate solution vectors, which are computed recursively at CG iterations. It is shown that when the number of CG iterations performed is less than the distance between two vertices in the graph, the round-off errors arising in the computation of one component do not affect the errors in the computation of another and are not correlated. At early iterations, errors propagate locally between neighboring vertices and are locally correlated. As the errors become sufficiently global, loss of orthogonality of residuals and conjugacy of search directions occurs, and stagnation of the iterative process sets in. The results of computational experiments are presented for matrices corresponding to graphs with a prescribed diameter, for Laplacians of differential equations, and for random matrices. These results demonstrate that the diameter and its approximate values, in particular the pseudodiameter, due to their low computational cost for sparse graphs, can be used in mixed-precision algorithms for solving linear systems, namely, for switching precision modes as round-off errors propagate and accumulate.

\end{abstract}

\section{Introduction}

The conjugate gradient method (hereinafter CG) is a fundamental iterative method for solving large systems of linear equations (hereinafter linear systems) $$Ax=b\eqno(1)$$ with a sparse symmetric positive definite matrix $A$. The CG method for solving the linear system (1) is constructed as a method of minimizing the quadratic form $$f(x)=\frac{1}{2} (Ax,x) - (b,x),$$ which attains its minimum at the solution of the system (1). If 
the computations are performed in exact arithmetic, that is, without round-off errors, then the CG method gives the solution of the linear system (1) with a matrix $A$ of size $n\!\times\!n$ in no more than $n$ iterations.

As a rule, the analysis of CG convergence, that is, estimation of its convergence rate, norms of error and residual vectors, etc., is performed based on the spectral characteristics of the matrix~$A$. In this approach, the matrix is considered as a linear operator. This approach is presented, for example, in~\cite{Carson}, where the convergence analysis of the method is reduced to considering linear systems with diagonal matrices, which is sufficient if only spectral characteristics are considered. The spectral characteristics of the matrix are directly related to the construction of CG as a method for solving linear systems: it is constructed as a method of finding a polynomial of the matrix~$A$ that attains its minimum on its spectrum. Minimization of the energy norm of the error vector is equivalent to finding a polynomial minimizing the weighted sum of squares of the polynomial values at the eigenvalues of the matrix. This principle underlies CG, which explains its efficiency for linear systems with matrices having a clustered spectrum. In this case, CG converges to the exact solution in a number of iterations close to the number of clusters.

Analyzing the metric properties of the graph of the matrix of a linear system, we consider it as the adjacency matrix of this graph, defining  connections between components of the vectors computed during CG iterations. The metric characteristics of the matrix graph are determined by distances between vertices of its graph. In this work, by the structure of a matrix we understand the structure of its weighted graph, whose vertices correspond to the row and column indices of the matrix, and edges --- to its nonzero entries. The matrix entries are considered by us as edge weights of the graph. The influence of the \textit{sparsity pattern} of a matrix (i.e., its unweighted graph) on the convergence of CG and other methods for solving the linear system (1) was previously considered in the contexts of preconditioning and the use of parallel computations for CG (for example, \cite{Saad}, \cite{Rude}).

We consider the impact of the structure of the weighted graph of the matrix on the accumulation of round-off errors in computations in finite-precision arithmetic. The basis of our approach is the analysis of the process of information propagation along paths in the matrix graph during the iterative process. \textit{Information propagation} occurs when updating the components of the vectors recursively computed at CG iterations. The vector components correspond to the vertices of the graph, and their influence on other components is determined by the edges of the graph. CG iterations are interpreted by us as iterations of traversals of the matrix graph, starting at the vertices corresponding to the nonzero components of the right-hand side of the linear system (the initial solution is assumed by us to be zero). We present a model of information and round-off error propagation during CG iterations, built on the basis of this approach. The paper shows that as long as the iteration number is less than the distance between two vertices in the graph, the information and round-off errors in the corresponding components do not affect each other.

By the \textit{correlation} of errors (in the sense of their structural coherence rather than statistical dependence), we understand two interrelated phenomena: \textit{spatial correlation}---the coherence of errors in different components of vectors at a single iteration, and \textit{temporal correlation}---the dependence of errors at later iterations on the errors that occurred at earlier ones. Both of these phenomena are caused by the deterministic propagation of round-off errors along the paths of the matrix graph during CG iterations.

At early CG iterations, information, including errors, propagates locally, that is, between nearby vertices, and the errors are locally correlated. This correlation is relatively small and has a minor effect on convergence. As iterations proceed, the correlation becomes increasingly global, that is, the error in computing a single component of a vector affects an increasing number of components. As a result, when the magnitude of these errors and their correlation become sufficiently large, loss of orthogonality of residuals and conjugacy of search direction vectors occurs.

Our proposed approach can be used, in particular, for mixed-precision computations organized as follows: computations critically important for achieving the required solution accuracy are performed in high precision, whereas the maximum possible fraction of computations is performed in the lower precision. The use of lower-precision floating-point formats reduces data transfer time and consumes fewer computational resources, which~provides a significant speedup of computations compared to computations performed exclusively in high precision. To illustrate this approach, we consider a two-stage mixed-precision computation scheme, in which at the first stage computations are performed using the CG method in single precision, and at the second---in double precision. By replacing a portion of the iterations performed in double precision with iterations performed in single precision mode, an approximate solution with the same accuracy as that obtained from computations performed entirely in double precision is achieved faster.

Hereinafter, by the \textit{computation precision mode} we mean the use of floating-point numbers of a certain precision. The metric characteristics of the matrix graph considered in this work---the diameter, the effective diameter, and the pseudodiameter---allow one to estimate a priori the moment 
of switching the computation precision mode. This moment is defined as the iteration number $N_{\mathrm{div}}$ starting from which the true residual $\|b\!-\! Ax^k\|$, computed explicitly, and the residual computed via the CG recurrence formula, differ in norm by at least an order of magnitude. This moment corresponds to the situation when switching to a mode with higher computation precision is necessary for further refinement of the approximate solution obtained up to the iteration preceding iteration $N_{\mathrm{div}}$. We put forward and confirm the following hypothesis: closeness of metric characteristics of graphs of matrices of the same type (for example, two Laplacians of partial differential equations obtained for different discretization parameters) corresponds to closeness of the values of $N_{\mathrm{div}}$ for them. One of the experiments presented below confirms this property as follows: for matrices of the same type with close diameters, the standard deviation of the value of $N_{\mathrm{div}}$ turns out to be small.

This work presents theoretical results on the correlation of round-off errors in the components of vectors recursively computed during CG iterations. The results of computational experiments are presented for matrices of the following types: adjacency matrices of unweighted graphs with introduced 
diagonal dominance, which most clearly reflect the influence of the matrix graph structure on CG convergence; Laplacians of partial differential equations; and random matrices with a prescribed condition number and diameter. These results show that the metric characteristics of the matrix graph 
allow one to detect and explain features of CG convergence that cannot be detected and explained solely by the spectral characteristics of the matrix.\\

\section{Spectral characteristics of the matrix \\and metric characteristics of the matrix graph}

\paragraph{Spectral characteristics of the matrix.} \textit{Spectral} characteristics of a matrix are the eigenvalues that make up the spectrum of the matrix, and the characteristics derived from them:
\begin{itemize}
\item[--] the set of \textit{eigenvalues} of the matrix $A$ --- the \textit{spectrum} of the matrix $A$, that is, the set of roots $\{\lambda_1,\ldots,\lambda_n\}$ of the \textit{characteristic equation} $$\det(A-\lambda I)=0,$$ where $I$ is the identity matrix;

\smallskip

\item[---] the number of distinct eigenvalues in the spectrum;

\smallskip

\item[---] the cluster structure of the spectrum;

\smallskip

\item[---] the \textit{condition number} of the matrix $\kappa(A)\!=\!\lambda_{\mathrm{max}}/\lambda_{\mathrm{min}}$, where $\lambda_{\mathrm{max}}$ and $\lambda_{\mathrm{min}}$ are the maximum and minimum eigenvalues from the spectrum of the matrix;

\smallskip

\item[---] the \textit{effective condition number} $\kappa_{\mathrm{eff}}(A)$, which determines the convergence rate of an iterative method taking into account the expansion of the right-hand side vector of the linear system in the eigenvector basis;

\smallskip

\item[---] the \textit{spectral radius} $\rho(T)\!=\!\max |\lambda_i(T)|$, where $T$ is the matrix obtained from the matrix~$A$ of the linear 
system for stationary iterative methods (Jacobi method, Gauss-Seidel method and other relaxation methods), and $\lambda_i(T)$ are the 
eigenvalues of $T$; 
\end{itemize}
and others.

The main advantage of spectral characteristics of the matrix from the point of view of CG convergence analysis is that they have a direct connection with the convergence process. By its construction, CG is a method that at the $k$-th iteration constructs a polynomial $p_k$ of degree not exceeding $k$, minimizing the weighted sum of squares of its values on the spectrum of the matrix. Namely, let $\hat{x}$ be the exact solution of the linear system (1), $x^k$ be the approximate solution at the $k$-th iteration (hereinafter, the symbol $\hat{\cdot}$ denotes the value of a quantity computed by CG 
in exact arithmetic), $e^k\!=\!\hat{x}\!-\!x^k$ be the error, $Ae^k\!=\!r^k$, and $\|e^k\|_A\!=\!\sqrt{(e^k,Ae^k)}$ be the energy norm of the error. Then at the $k$-th iteration, we solve the following optimization problem: $$\|e^k\|_A^2=\sum\limits_{i=1}^nw_i[p_k(\lambda_i)]^2\to\min,$$ where the minimum is sought subject to the constraint $p_k(0)\!=\!1$, where $w_i\!=\!\lambda_i c_i^2$, and $c_i$ are the coordinates of the initial error vector $e^0\!=\!\hat{x}\!-\!x^0$ in the eigenvector basis of the matrix $A$. Although this polynomial is not constructed explicitly during CG iterations, 
it underlies the theoretical analysis of its convergence. Therefore, spectral characteristics allow one to obtain theoretical estimates of the CG convergence rate, the number of CG iterations required, and the accuracy of the resulting approximate solution, including estimates of the norms of the 
error vector and the residual vector.

For the dynamic choice of the computation precision mode during the iterative solution of linear systems, the practical use of spectral characteristics faces the following difficulties. Computing (estimating) the spectral characteristics of the matrix cannot be performed with sufficient accuracy in a time negligible compared to the running time of CG itself. Computing the entire spectrum has complexity $O(n^3)$ and fundamentally cannot be used 
to employ information about the cluster structure of the spectrum for choosing the precision mode. Estimates of the maximum and minimum eigenvalues obtained at early iterations via Ritz values at early iterations turn out to be unreliable: the smallest Ritz value systematically overestimates $\lambda_{\mathrm{min}}$, especially if there are clusters close to zero. This leads to an underestimated condition number and an unjustified switch to a 
lower-precision computation mode, or to a delayed switch to a higher-precision computation mode.

Spectral characteristics are \textit{global}: they are determined by the entire matrix and only implicitly reflect its \textit{local} characteristics, that is, the collection of values of elements in individual rows (columns) and the dependencies between them defined by the matrix of the linear  system. This global nature of spectral characteristics makes them difficult to interpret when considering a particular CG convergence process: two structurally different linear system matrices with the same spectrum are indistinguishable under this approach.\\

\paragraph{Metric characteristics of the matrix graph.} \textit{Metric characteristics of the matrix graph} of an $n\!\times\!n$ matrix $A$ are the characteristics determined by the graph whose adjacency structure is defined by the nonzero entries of the matrix $A$. The \textit{graph of the matrix} $A$ is defined as a graph $G\!=\!G(A)\!=\!(V,E)$ with a set of vertices $V\!=\!\{1,\ldots,n\}$, where vertex $i\!\in\!V$ is associated with the $i$-th row (column) of the matrix, and the set of edges $E\!=\!\{(i,j)\, |\, a_{ij}\!\neq\! 0\}$ corresponds to the positions of the nonzero elements of the matrix (including self-loops corresponding to the diagonal elements). Let $m\!=\!|E|$ denote the number of edges in the matrix graph. Such a graph is referred to as the \textit{sparsity pattern} of the matrix.

A sparse matrix corresponds to a sparse graph, that is, one in which $m\!=\!O(n)$. The graph of a sparse matrix can be considered as \textit{unweighted}, in which all edge weights are equal to $1$, or as \textit{weighted}, where the weight of an edge $(i,j)\!\in\!E$ is equal to the matrix element 
$a_{ij}\!\in\!\mathbb{R}$. In what follows, we consider weighted matrix graphs.

In addition, the graphs considered below are \textit{connected}, that is, any vertex $i\!\in\!V$ is connected to any other vertex $j\!\in\!V$ by some sequence of edges $\{(i_0,i_1),(i_1,i_2)\dots,(i_{k-1},i_k)\}$, where $(i_{t-1},i_t)\!\in\!E$, $t\!=\!1,\dots,k$, $i_0\!=\!i$, $i_k\!=\!j$. This sequence defines a \textit{walk of length $k$} in the graph $G$. By a \textit{walk in a graph}, we hereafter mean a sequence of vertices (or incident edges) in which vertices (edges) may repeat. The length of a walk is equal to the number of edges in it, including self-loops incident to vertices. 
Let $\ell(i,j)$ denote the \textit{distance} between vertices $i$ and $j$ in the graph, which is defined as the length of the shortest walk connecting them. Information and round-off errors propagate along walks in the graph.

Matrices with disconnected graphs can be reduced to a block-diagonal form, where each block and its corresponding graph can be considered separately, in the same way as will be done below.

\smallskip

The metric characteristics considered in this work are as follows:

\begin{itemize}
\item[---] the distance between vertices in the graph;

\smallskip

\item[---] the \textit{diameter} $\mathcal{D}$ of the matrix graph: $$\mathcal{D}=\max\limits_{i,j\in V}\{\ell(i,j)\};$$

\smallskip

\item[---] the \textit{effective diameter} $\mathcal{D}_{\mathrm{eff}}$ --- the distance at which the majority of pairs of vertices in the matrix graph are separated; for a given $p\!\in\!(0,1)$, $$\mathcal{D}_{\mathrm{eff}}\!=\!\min\biggl \{l\in\mathbb{N}\, \biggl |\, \frac{|\{(i,j)\, :\, \ell(i,j)\le l\}|}{n(n-1)}\ge p\biggr \};$$ for example, for $p\!=\!0.9$ this is the smallest value of $l$ such that at least $90\%$ of all pairs of vertices are at a distance not exceeding $l$ (this value of $p$ is used by us hereafter in computations of $\mathcal{D}_{\mathrm{eff}}$);

\smallskip

\item[---] the \textit{pseudodiameter} $\tilde{\mathcal{D}}$ of the matrix graph --- a two-stage approximation of the diameter $\mathcal{D}$: at the first stage, a vertex $j$ is found that is maximally distant from a randomly chosen vertex $i$; at the second stage, a vertex $t$ is found that is maximally distant from $j$, and the value of the pseudodiameter $\tilde{\mathcal{D}}$ is set to $\ell(j,t)$; this two-stage process can be repeated 
several times to refine the value of~$\tilde{\mathcal{D}}$ as an estimate of $\mathcal{D}$;

\end{itemize}

\smallskip

The complexity of computing the diameter $\mathcal{D}$ is $O(n^2)$, whereas its approximate value and the pseudodiameter can be computed for sparse matrices in subquadratic time in $n$. A $3/2$-approximation of the diameter can be computed in time $O(m\sqrt{n})$ \cite{RodittyVassillevska}, which for sparse graphs amounts to $O(n^{3/2})$. The pseudodiameter $\tilde{\mathcal{D}}$ is one of the most rapidly computable estimates of the diameter, which can also serve as an upper bound for the effective diameter. Computing the distance between vertices in the graph and, accordingly, computing the 
pseudodiameter, can be performed using breadth-first search (BFS). The pseudodiameter $\tilde{\mathcal{D}}$ can be regarded as an approximation of the diameter $\mathcal{D}$. Finding it requires time not exceeding the time needed to perform several CG iterations \cite{Prolubnikov2}, and for sparse graphs amounts to $O(m)$, since the computational complexity of a single BFS traversal for sparse graphs is $O(m)$ \cite{Burkhardt}.

An estimate of the effective diameter $\mathcal{D}_{\mathrm{eff}}$ can be obtained in $O(s\cdot \mathcal{D}\cdot m)$ time \cite{Amati}, where $s$ is the number of randomly selected vertices used to refine the estimate. 

The metric characteristics of the matrix graph are intuitive and admit a natural interpretation. They characterize the structure of the matrix and~the interaction between the components of the vectors computed at the iterations, describing the \textit{local} structure of these dependencies. Their fundamental influence on CG convergence is scaled by the numerical values of the matrix entries (the edge weights of the graph), which determine the weights of the walks along which both information and round-off errors propagate.

The main disadvantage of such matrix characteristics for CG convergence analysis is that their connection with CG is indirect, since they are not used directly at CG iterations. This makes it difficult to obtain accurate theoretical estimates of the norms of error and residual vectors, the number of iterations required for convergence.

\smallskip

For arbitrary sparse matrices, the spectral and metric characteristics of the matrix graph are related implicitly and ambiguously, in contrast to adjacency matrices and Laplacians of graphs, for which this connection is rather well established and studied. A large number of works are 
devoted to this topic, among which the monograph \cite{Cvetkovic} can be noted as the most widely known and comprehensive one.\\

\section{Conjugate gradient method\\ as a matrix graph traversal algorithm}

\subsection{Metric characteristics of the matrix graph and CG convergence}

\textit{Breadth-first search (BFS)} is an algorithm for traversing an undirected graph that visits vertices layer by layer in increasing order of distance from the starting vertex. BFS is guaranteed to find the shortest paths. In \cite{Prolubnikov1}, it is shown that the Jacobi method, as well as the Gauss-Seidel method and other relaxation methods for solving linear systems, can be considered as matrix graph traversal algorithms. One iteration of the Jacobi method corresponds to one iteration of BFS, whereas relaxation methods, in particular the Gauss-Seidel method, induce traversals different from BFS because the propagation of information---namely, the values of the components of the approximate solution---occurs differently: unlike the Jacobi method, at each $k$-th iteration, the components already computed during this same iteration are taken into 
account.

The main result of \cite{Prolubnikov1} is the following Theorem~1, which allows one to consider these linear system solvers as graph traversals. Consider the matrix $A$ of the linear system (1) as the adjacency matrix of a weighted undirected graph $G\!=\!G(A)$. Let $s\!\in\! V(G)$ be the starting vertex of the traversal, $x^0\!=\!(0,\dots,0)$ be the initial approximation, and the right-hand side $b$ of the linear system be the vector having a single nonzero component in the $s$-th position, that is, $b\!=\!e_s$, where $e_s$ is the $s$-th vector of the standard basis in $\mathbb{R}^n$. Let $x^k$ denote the approximate solution at the $k$-th CG iteration. (Note that hereinafter $e$ does not denote the error vector.) The following holds

\bigskip

\noindent\textbf{Theorem 1} \cite{Prolubnikov1}\textbf{.} \textit{The condition $x_j^{k-1}\!=\!0$ and $x_j^k\!\neq\!0$ holds if and only if vertex $j\!\in\! V$ is reached for the first time in the traversal of the graph $G$ from vertex $s$ at the $k$-th iteration of the traversal, corresponding to 
the $k$-th iteration of the numerical method.}

\bigskip

Consider the scheme of the CG iterative process. Let $N$ 
denote the number of iterations performed.

\begin{table}[h!]

\begin{tabular}{lll}
0. & $k\!=\!0$; $d^{0}\!=\!r^{0}\!=\!b\!-\!Ax^{0}$; & --- \textit{initialization}\\
\quad & \textbf{for} $k\!:=1$ \textbf{to} $N$\\
1. &\qquad $\alpha_{k}\!=\!\bigl (r^{k},r^{k}\bigr )/\bigl (d^{k}, Ad^{k}\bigr )$; & --- \textit{step size} \\
2. &\qquad $x^{k+1}\!=\!x^{k}+\alpha_{k}d^{k}$; & --- \textit{approximate solution}\\
3. &\qquad $r^{k+1}\!=\!r^{k}\!-\!\alpha_{k}Ad^{k}$; & --- \textit{recursively computed residual}\\
4. &\qquad $\beta_{k+1}\!=\!\bigl (r^{k+1},r^{k+1}\bigr )/\bigl (r^{k},r^{k}\bigr )$; & \\
5. &\qquad $d^{k+1}\!=\!r^{k+1}+\beta_{k+1}d^{k}$ & --- \textit{direction (search direction)} \\
\end{tabular}
\end{table}

\noindent We assume $x^0\!=\!(0,\dots,0)\in\mathbb{R}^n$. For better readability, we use subscripts for the scalar coefficients $\alpha_k,\beta_k\!\in\!\mathbb{R}$.

Below, we show that CG for a linear system with the right-hand side $b\!=\!e_s$ can also be interpreted as a traversal of the matrix graph starting at vertex $s$. For an arbitrary right-hand side $b$, the value of a component of the search direction vector $d^k$ can be represented as a sum of \textit{contributions} --- values coming from the nonzero components of $b$: each contribution corresponds to a walk traversed during a graph traversal starting at a vertex $i\!\in\!V$ such that $b_i\!\neq\!0$, and each component of $d^k$ is the sum of contributions arriving at it along walks 
ending at it. As we show below, the nonzero components in $d^k$ can correspond only to those vertices that are reachable in no more than $k$ iterations of a traversal starting at vertex $i\!\in\!\mathrm{supp}(b)$, where $\mathrm{supp}(\cdot)$ denotes the set of indices of nonzero components of a vector.

For the convergence of the numerical implementation of CG, the following two equivalent conditions, which hold in exact arithmetic, must be satisfied with sufficient accuracy at each iteration:
\begin{itemize}
\item[(a)] the recursively computed residuals $r^k$ are mutually orthogonal: $$(r^{k_1},r^{k_2})\!=\!0,\quad k_1\!\neq\!k_2;$$
\item[(b)] the search direction vectors $d^k$ are mutually conjugate \textit{($A$-orthogonal)}: $$(d^{k_1}, Ad^{k_2})\!=\!0,\quad  k_1\neq k_2.$$
\end{itemize}

Due to the limited length of the mantissa of machine numbers, round-off errors appear during the iterations. Through the CG recurrence relations, these errors propagate between the components of the residual and search direction vectors. The increase in the magnitude and correlation of these errors leads to a loss of orthogonality of the residuals. As a result, convergence slows down, followed by \textit{stagnation} (\textit{stagnation of the iterative process} or \textit{residual stagnation}): the norms of the error $\hat{x}\!-\!x^k$ and the \textit{true residual} $r_{\mathrm{true}}^k\!=\!b-Ax^k$ cease to decrease. Let $r_{\mathrm{rec}}^k$ denote the \textit{recursive residual}, i.e., the vector $r^k$ computed at step 3 of the CG scheme 
above. In exact arithmetic, $r_{\mathrm{true}}^k\!=\!r_{\mathrm{rec}}^k$. In a numerical implementation, starting from a certain iteration, $r_{\mathrm{true}}^k$ and $r_{\mathrm{rec}}^k$ begin to differ increasingly from each other.

\begin{figure}[htbp]
\begin{center}
\includegraphics[width=85mm]{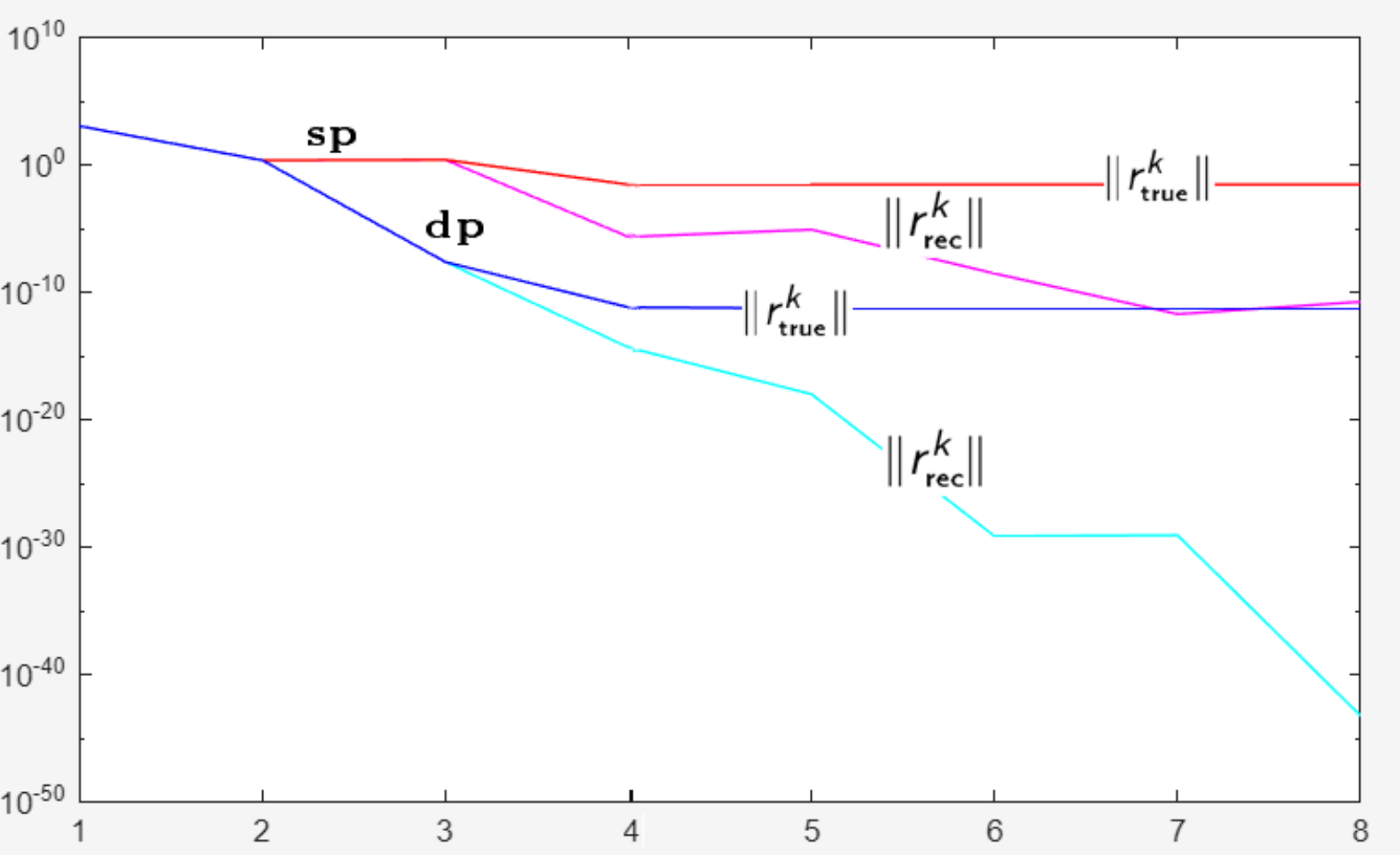}
\caption{CG convergence for a matrix with a star graph, $n\!=\!600$.}\label{fig:1}
\end{center}
\end{figure}

Figure~\ref{fig:1} shows, on a semi-logarithmic scale, the dependence of $\|r^k\|$ ($Y$-axis) on the iteration number $k$ ($X$-axis) when solving a linear system in single precision ($\mathrm{\textit{sp}}$) and double precision ($\mathrm{\textit{dp}}$) for a $601\!\times\!601$ matrix with a star graph structure. This graph consists of a single central vertex connected to $600$ pairwise non-adjacent vertices; thus, the values of the 
considered metric characteristics are minimal for this graph among all graphs of this dimension. As the plots show, for both precision modes, a significant difference between $\|r_{\mathrm{true}}^k\|$ and $\|r_{\mathrm{rec}}^k\|$ is observed as early as after the 3rd iteration. In subsequent 
iterations, the accumulation of errors and their correlation lead to a loss of orthogonality of the residual vectors and to stagnation.

Note that the anomalously small values of the norms of the recursively updated residual, reaching $10^{-10}$ in single precision and unacceptably small values (down to $10^{-40}$ and below) in double precision, are caused by a catastrophic loss of significant digits. For the star graph, which has 
only three distinct eigenvalues, CG in exact arithmetic converges in no more than 3 iterations. In subsequent iterations, the true residual $\|r_{\mathrm{true}}^k\|$ ceases to decrease, but the recursively computed residual continues to do so: the vectors $r_{\mathrm{rec}}^k$ and 
$\alpha_kAd^k$ become nearly linearly dependent, as do their errors. Subtracting the latter from the former in floating-point arithmetic using denormalized numbers ($10^{-20}, 10^{-30}, \dots$) causes the recursive residual to continue decreasing, whereas in reality this only reflects 
accumulated computation errors. This effect is amplified by the fact that, when performing vector-matrix operations in the MATLAB environment, intermediate results are not rounded to the precision of the operand format, which reduces the level of round-off errors introduced at each iteration.

\subsection{Are spectral characteristics of the matrix \\ sufficient to characterize CG convergence?}

That spectral characteristics alone are insufficient to fully characterize CG convergence for a given matrix and a particular right-hand side can be illustrated by the following examples. Consider the CG convergence for pairs of matrices consisting of a non-diagonal matrix $A$ and a diagonal matrix 
$D(A)$, whose nonzero entries are the eigenvalues of $A$. For a matrix $A$ with spectrum $\{\lambda_1,\dots,\lambda_n\}$, we solve the linear system $Ax\!=\!b$ with a known exact solution $\hat{x}\!=\!(1,\dots,1)$ and right-hand side $b\!=\!A\hat{x}$, and also solve the linear system $D(A)x\!=\!b$, where $D(A)\!=\!\mathrm{diag}(\lambda_1,\dots,\lambda_n)$.

\begin{figure}[htbp]
\begin{center}
\includegraphics[width=125mm]{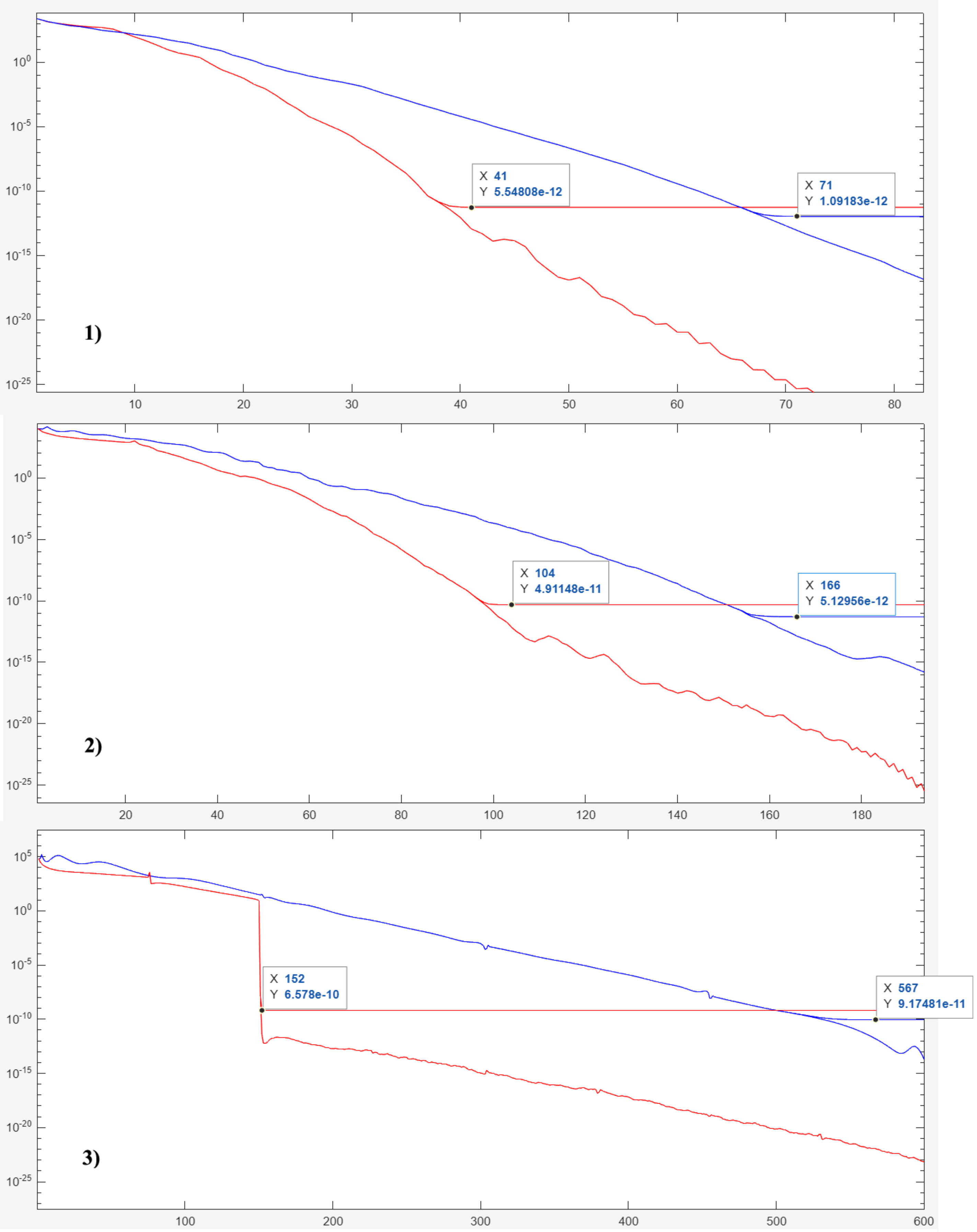}
\caption{CG convergence for Laplacian matrices and diagonal matrices with the same spectrum.}\label{fig:2}
\end{center}
\end{figure}

\begin{figure}[htbp]
\begin{center}
\includegraphics[width=125mm]{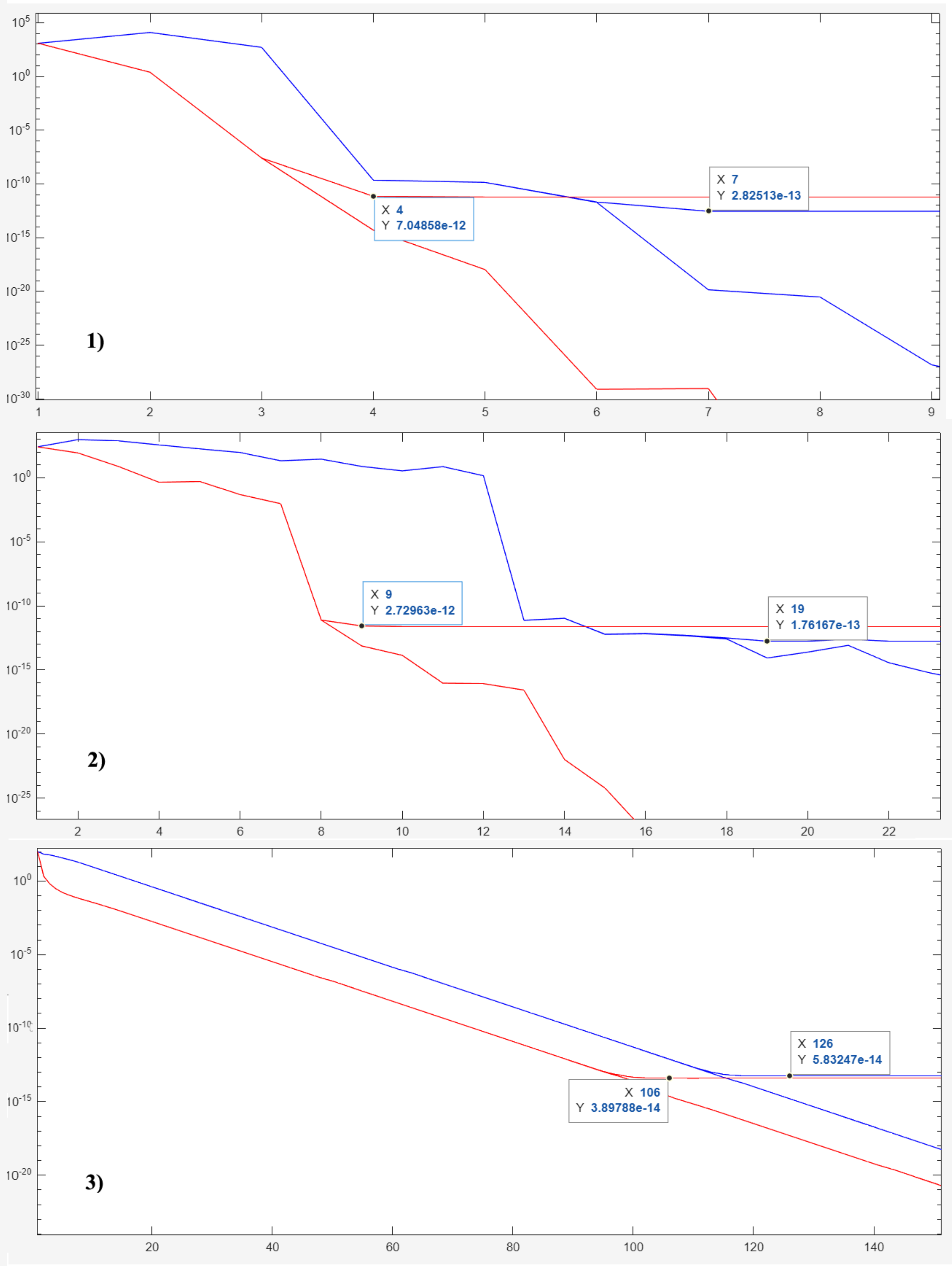}
\caption{CG convergence for matrices with extended star graphs and diagonal matrices with the same spectrum.}\label{fig:3}
\end{center}
\end{figure}

Figures~\ref{fig:2} and~\ref{fig:3} show the plots of the norms of the true residual $r_{\mathrm{true}}^k\!=\!b-Ax^k$ and the recursively computed residual $r_{\mathrm{rec}}^k$ ($r_{\mathrm{rec}}^k\!=\!r_{\mathrm{rec}}^{k-1}\!-\!\alpha_kAd^k$) versus the iteration number. The plots for non-diagonal matrices are shown in red (their plots appear further to the left in all figures), and those for the corresponding diagonal matrices with eigenvalues on the diagonal --- in blue. They show that the CG convergence for two such matrices differs significantly in its behavior.

A diagonal matrix corresponds to a disconnected graph consisting of a set of isolated vertices with weighted self-loops. The diameter of a disconnected graph is considered to be infinite. The smaller the diameter (effective diameter) of the matrix graph, the greater the possible differences in CG convergence for linear systems with matrices $A$ and those with matrices $D(A)$. This is caused by the fact that the smaller the distances between 
vertices in the matrix graph, the faster round-off errors accumulate, and the more correlated the errors become at different iterations, which leads to a faster loss of conjugacy of directions and orthogonality of residuals, and, as a result, to stagnation of convergence.

In the case of diagonal matrices, the accumulation of errors during CG occurs only in the computation of inner products and componentwise arithmetic operations. This fundamentally distinguishes the accumulation of round-off errors during CG iterations for diagonal matrices from the accumulation of 
errors for non-diagonal ones. In this case, there is no propagation of information between the components of the vectors computed at the iterations: the product $Ad^k$ of the direction vector by the diagonal matrix $A$ only scales the components of the vector, but the values of the components of $d^k$ and the round-off errors obtained in their computation do not affect each other. The same holds for the approximate solution and residual vectors, as well as for the corresponding error vectors. As a result, round-off errors do not propagate between components, so they are uncorrelated, and stagnation of the iterative process occurs later.

We say that \textit{residual divergence} occurs at some iteration $k$ if $$\frac{\|r_{\mathrm{true}}^k\|}{\|r_{\mathrm{rec}}^k\|}>10,$$ that is, the norm of the true residual exceeds the norm of the recursively computed residual by a factor of 10. We denote this iteration by $N_{\mathrm{div}}$. Here, the value 10 serves as a threshold separating the convergence phase from the stagnation phase, at which point, due to accumulated round-off errors, the ratio $\|r_{\mathrm{true}}^k\|/\|r_{\mathrm{rec}}^k\|$ begins to grow exponentially.

Figures~\ref{fig:2} and~\ref{fig:3} show the plots of the true and recursively computed residuals versus the CG iteration for the matrices $A$ and $D(A)$ computed in double precision. The iteration numbers $N_{\mathrm{div}}$ are indicated. Figure~2 shows the plots for Laplacian matrices of differential equations, $n\!=\!600$; on the right --- for matrices with extended star graphs, $n\!=\!601$. The parameters of the matrices $A$ are as follows:

\begin{itemize}
\item[--] Fig.~2.1: $\mathcal{D}\!=\!23$, $\tilde{\mathcal{D}}\!=\!23$, $\mathcal{D}_{\mathrm{eff}}\!=\!13$;
\item[--] Fig.~2.2: $\mathcal{D}\!=\!53$, $\tilde{\mathcal{D}}\!=\!53$, $\mathcal{D}_{\mathrm{eff}}\!=\!33$;
\item[--] Fig.~2.3: $\mathcal{D}\!=\!152$, $\tilde{\mathcal{D}}\!=\!152$, $\mathcal{D}_{\mathrm{eff}}\!=\!104$;

\smallskip

\item[--] Fig.~3.1: $\mathcal{D}\!=\!2$, $\tilde{\mathcal{D}}\!=\!2$, $\mathcal{D}_{\mathrm{eff}}\!=\!2$;
\item[--] Fig.~3.2: $\mathcal{D}\!=\!10$, $\tilde{\mathcal{D}}\!=\!10$, $\mathcal{D}_{\mathrm{eff}}\!=\!9$;
\item[--] Fig.~3.3: $\mathcal{D}\!=\!600$, $\tilde{\mathcal{D}}\!=\!600$, $\mathcal{D}_{\mathrm{eff}}\!=\!411$.
\end{itemize}

\noindent As can be seen, in all six cases CG converges faster for non-diagonal matrices, and as the values of $\mathcal{D}$, $\tilde{\mathcal{D}}$, $\mathcal{D}_{\mathrm{eff}}$ increase for matrices of both types, stagnation occurs later.\\

\section{Propagation of information along the matrix graph \\ during conjugate gradient iterations}

The iterative solution of linear systems by the conjugate gradient method can be viewed as the propagation of \textit{information} along the matrix graph---the values of the vector components computed at the iterations of the algorithm, that is, in the case of CG, the components of the 
approximate solution, residual, and search direction. The propagated information includes \textit{errors} caused by computing these values with round-off errors in finite-precision arithmetic. Hereinafter, we consider the expressions \textit{``propagates between vector components 
with indices $i$ and $j$''} and \textit{``propagates between vertices $i$ and $j$ of the matrix graph''} to be equivalent.\\

\subsection{Propagation of information along the matrix graph \\ during its breadth-first traversal}

Consider the following iterative process --- \textit{algebraic BFS}, which implements a breadth-first traversal of the graph $G$ as a sequence of vector-matrix multiplications by the adjacency matrix $A$ of the graph: $$b^k=Ab^{k-1}, \eqno(2)$$ where, for the starting vertex $i$ of the traversal, $b^0$ is the vector with a single nonzero component $b_i^0$. Vertex $j$ is reachable from vertex $i$ at the $k$-th BFS iteration if $b_j^{k-1}\!=\!0$ and $b_j^k\!\neq\!0$.

\begin{figure}[htbp]
\begin{center}
\includegraphics[width=27mm]{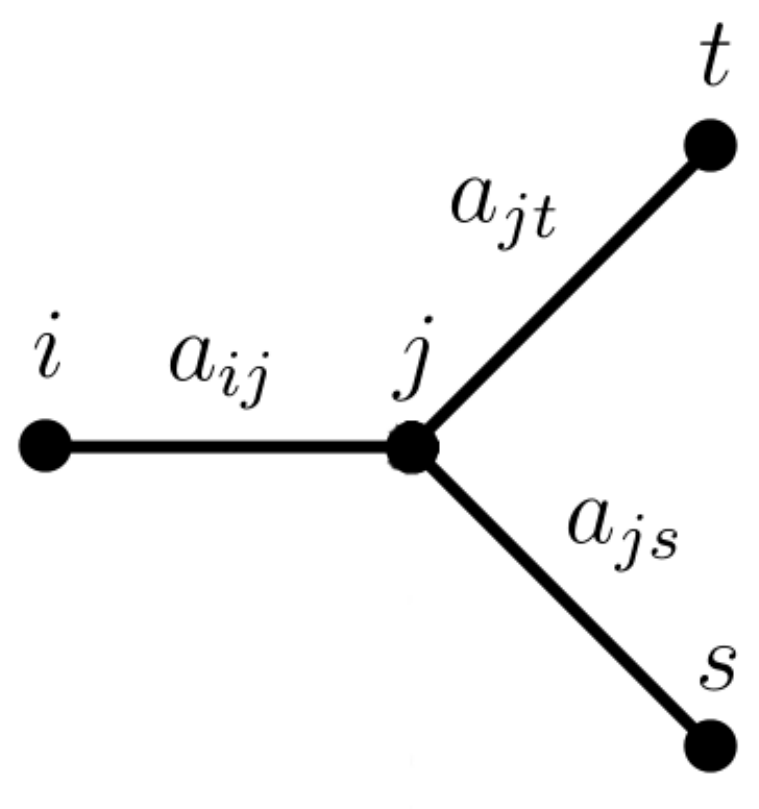}
\caption{A subgraph on 4 vertices of a certain graph.}\label{fig:4}
\label{fig1}
\end{center}
\end{figure}

For a subgraph on 4 vertices of a certain graph in Fig.~\ref{fig:4}, in which vertex $j$ is adjacent to vertices $i$, $t$, $s$, we have: $$b_j^1=\sum\limits_{q=1}a_{jq}b_q^0=\sum\limits_{(j,q)\in E}a_{jq}b_q^0=a_{ji}b_i^0+a_{jt}b_j^0+a_{js}b_s^0=a_{ji}b_i^0,$$ since $b_t^0\!=\!b_s^0\!=\!0$ and $a_{ji}\!\neq\!0$, $b_i^0\!=\!1$. Consequently, $b_j^1\!\neq\!0$, whereas $b_t^1\!=\!b_s^1\!=\!0$, since $b_t^1\!=\!a_{jt}b_j^0\!=\!0$ and $b_s^1\!=\!a_{js}b_s^0\!=\!0$. Therefore, after the first iteration of the algebraic BFS, the corresponding components of the vector $b^1$ will be nonzero only for vertices adjacent to vertex $i$. Such a vertex is $j$, whereas the $t$-th and $s$-th components of $b^1$ remain zero. They will become 
nonzero only after the second iteration, since $\ell(i,t)\!=\!\ell(i,s)\!=\!2$: $b_t^2\!=\!a_{jt}b_j^1\!\neq\!0$, $b_s^2\!=\!a_{js}b_j^1\!\neq\!0$.

Continuing the iterations of the algebraic BFS for subsequent values of $k$, we obtain $$b_j^k=\sum\limits_{w\in W_{ij}^k}C_{\mathrm{BFS}}(w),\eqno(3)$$ 
where $W_{ij}^k$ is the set of walks $w$ of length $\ell(w)\!\le\! k$ starting at vertex $i$ and ending at vertex $j$, and $C_{\mathrm{BFS}}(w)$ is the \textit{contribution} delivered to vertex $j$ along the walk $w\!=\!\{(i_0,i_1),\dots (i_{\ell(w)-1},i_{\ell(w)})\}$, $i_0\!=\! i$, $i_{\ell(w)}\!=\!j$: 
$$C_{\mathrm{BFS}}(w)=(a_{i_0i_1}a_{i_1i_2}\cdots a_{i_{\ell(w)-1}i_{\ell(w)}})b_i^0=\biggl (\prod\limits_{t=1}^{\ell(w)}a_{i_{t-1}i_t}\biggr )b_i^0.$$

\subsection{Propagation of information along the matrix graph \\ 
during conjugate gradient iterations}

CG iterations can also be considered as iterations of a traversal of the matrix graph. Let $\mathrm{supp}(b)\!=\!\{i\}$. Since for $k\!>\!1$
$$r^k=r^{k-1}-\alpha_kAd^{k-1},\quad d^k=r^k+\beta_k d^{k-1},$$ $$r^{k-1}=d^{k-1}-\beta_{k-1}d^{k-2},$$ we have $$d^k=(1+\beta_k)d^{k-1}-\beta_{k-1}d^{k-2}-\alpha_{k-1}Ad^{k-1},\eqno(4)$$ which in componentwise notation takes the form $$d_j^k=(1+\beta_k)d_j^{k-1}-\beta_{k-1}d_j^{k-2}-\alpha_{k-1}\sum\limits_{t=1}^n a_{jt}d_t^{k-1}.\eqno(5)$$ The summation in (5) of products with nonzero values of $a_{jt}$ is equivalent to summing the products $a_{jt}d_t^{k-1}$ over all vertices $t$ adjacent to $j$, that is, $$d_j^k=(1+\beta_k)d_j^{k-1}-\beta_{k-1}d_j^{k-2}-\alpha_{k-1}\sum\limits_{(j,t)\in 
E}a_{jt}d_t^{k-1}.\eqno(6)$$ The notation ``$(j,t)\!\in\!E$'' is equivalent to the notation ``$\{t\, |\, (j,t)\!\in\! E\}$''.

According to (6), the traversal of the matrix graph during CG is realized through the term $-\alpha_{k-1}\sum_{(j,t)\in E}a_{jt}d_t^{k-1}$, where according to the CG formulas $\alpha_{k-1}\!\neq\!0$. Before reaching the iteration with number $k\!=\!\ell(i,j)$, we have $d_j^{k-2}\!=\!0$ and $d_j^{k-1}\!=\!0$. Consequently, before this iteration, (6) takes the form: $$d_j^k=-\alpha_{k-1}\sum\limits_{(j,t)\in E}a_{jt}d_t^{k-1},$$ where $d^0\!=\!b$. Similarly to BFS, this means that while $k\!<\!\ell(i,j)$, there are no traversed walks in the matrix graph for propagating information---the value $d_i^0$---to vertex $j$ and contributing to the value $d_j^k$. Consequently, $d_j^k$ can be nonzero only if $\ell(i,j)\!\le\!k$. We formalize this observation later.

\bigskip

\noindent \textbf{Remark 1.} Unlike the stationary iterative methods considered in~\cite{Prolubnikov1}, for which, by Theorem~1, when $\mathrm{supp}(b)\!=\!\{i\}$, it is guaranteed that the component of the approximate solution vector is nonzero upon its first reach along any walk, for CG it is 
theoretically possible that $d_j^k\!=\!0$ for $k\!\ge\!\ell(i,j)$. This can occur in exact arithmetic because contributions from vertex~$i$ to vertex~$j$ along different walks may have different signs and sum to zero. However, in finite-precision arithmetic, such theoretically possible mutual cancellation of contributions is unlikely, since the errors at a CG iteration are not only the errors propagated from earlier iterations, but also new round-off errors arising from computations involving the nonzero contribution values obtained at the current iteration. This error is of order $\mathcal{O}(\varepsilon_{\mathrm{mach}})$, where $\varepsilon_{\mathrm{mach}}$ is the machine epsilon. Therefore, achieving exact equality $d_j^k\!=\!0$ is  practically impossible.\\

\subsection{Information propagation model with additive error}

Let us denote the error that we have in the search direction vector $d^k$ after performing $k$ CG iterations by $\epsilon^k$: $$d^k=\hat{d}^k+\epsilon^k.$$ The error $\epsilon^k$ can be represented as $$\epsilon^k=\sum\limits_{t=1}^k\varepsilon^{t,k},\eqno(7)$$ where $\varepsilon^{t,k}$ is the error obtained by propagating the error $\varepsilon^t\!=\!\varepsilon^{t,t}$, which includes:
\begin{itemize}
\item[(a)] the error arising from performing arithmetic operations with round-off errors at the $t$-th iteration;
\item[(b)] the error arising from the propagation along the matrix graph of errors that occurred at previous iterations. 
\end{itemize}

\noindent In particular, at the first iteration $\epsilon^1\!=\!\varepsilon^1=\!\varepsilon^{1,1}$ and $$d^1=\hat{d}^1+\varepsilon^1.$$ At each subsequent $k$-th iteration, $\varepsilon^1$ will produce the error $\varepsilon^{1,k}$ contributing to~(7).

Using the additive error model, we replace all computation errors at the $k$-th CG iteration by an error vector that is added to the search direction vector, assuming that the values $\alpha_k$ and $\beta_{k+1}$ are computed exactly. If $\delta\alpha_k$ is the error in computing the exact value 
$\hat{\alpha}_k$, that is, $\alpha_k\!=\!\hat{\alpha}_k+\delta\alpha_k$, then according to (4), $\delta\alpha_k$ affects~$d^{k+1}$ through the term $\delta\alpha_k\cdot Ad^k$. This means that the error $\delta\alpha_k$, obtained in computing $\alpha_k$, propagates to the components of $d^{k+1}$ corresponding to the nonzero components of the vector $Ad^k$. Similarly, errors in computing $\beta_k$ and $\beta_{k-1}$ will not lead to the 
appearance of new nonzero components in $d^{k+1}$, as follows from (4). This means that errors in~$\alpha_k$ and~$\beta_{k+1}$ do not introduce new dependencies between round-off errors at CG iterations, but only scale them.

Thus, it can be asserted that there exists a sequence $\epsilon^k$, constructed from the sequences $\varepsilon^{t,k}$, such that CG, in which the values $\alpha_t$ and $\beta_{t+1}$ are assumed to be computed exactly, with their computation errors accounted for in the additive error, produces the same values $x^t$, $r^t$, $d^{t+1}$ at the $t$-th iteration as the computations with errors in~$\alpha_t$ and~$\beta_{t+1}$. Consequently, the 
propagation of errors obtained in computing inner products at CG iterations will occur through the vector $\varepsilon^t$, which includes them as an additive component.

\smallskip

Consider the propagation of the round-off error present only after the first iteration---the error $\varepsilon^1$---along the matrix graph during CG iterations. The propagation of errors $\varepsilon^t$, $t\!>\!1$, proceeds similarly.

We have $$d_j^2=(1+\beta_2)d_j^1-\beta_1d_j^0-\alpha_1\sum_{(j,t)\in E}a_{jt}d_t^1=$$ $$=(1+\beta_2)(\hat{d}_j^1+\varepsilon_j^1)-\beta_1(d_j^0+\varepsilon_j^0)-\alpha_1\sum_{(j,t)\in E}a_{jt}(\hat{d}_t^1+\varepsilon_t^1)=$$ $$=\biggl ((1+\beta_2)\hat{d}_j^1-\beta_1d_j^0-\alpha_1\sum_{(j,t)\in 
E}a_{jt}\hat{d}_t^1\biggr )+$$ $$+\biggl ((1+\beta_2)\varepsilon_j^1-\beta_1\varepsilon_j^0-\alpha_1\sum_{(j,t)\in E}a_{jt}\varepsilon_t^1\biggr ),\eqno(8)$$ where $\varepsilon_j^0\!=\!0$ if the assignment operations during initialization (step 0 of CG) are performed without loss of precision. The loss is possible if the assignment operations are executed in a lower precision than the precision of the data fed to the input (the matrix and the right-hand side of the linear system).

Thus, from (8) we have $$d_j^2=\hat{d}_j^2+\varepsilon_j^2,$$ where $\hat{d}_j^2$ is the exact value of the $j$-th component of the search direction vector after the second iteration, and the round-off error obtained at the first iteration propagates to the second one as $\varepsilon_j^2$: $$\varepsilon_j^2=(1+\beta_2)\varepsilon_j^1-\beta_1\varepsilon_j^0-\alpha_1\sum_{(j,t)\in E}a_{jt}\varepsilon_t^1.$$ 

Continuing analogously, due to the linearity of the CG transformation (6), we obtain the formula for $\varepsilon_j^k$---the error arising from the propagation of $\varepsilon^1$ along the matrix graph to vertex $j$ after $k$ iterations: $$\varepsilon_j^k=(1+\beta_k)\varepsilon_j^{k-1}-\beta_{k-1}\varepsilon_j^{k-2}-\alpha_{k-1}\sum\limits_{(j,t)\in E}a_{jt}\varepsilon_t^{k-1},\eqno(9)$$ where $\varepsilon_t^{k-1}\!=\!\varepsilon_t^{1,k-1}$.

\smallskip

Note that above, in deriving (9), we considered the propagation during $k$ iterations of only a single error $\varepsilon^1$---the error that arises at the first iteration. Similarly, the error $\varepsilon^t$, arising at the $t$-th iteration, $t\!>\!1$, at a graph vertex (vector 
component) as a result of performing arithmetic operations at the current iteration, propagates further along the same walks in the graph that were traversed during the propagation of the error $\varepsilon^1$. Due to the linearity of transformation~(9) with respect to the error vectors, the 
total error at the $k$-th iteration is the sum (7) of all such error vectors propagated along the graph at iterations up to the $k$-th and arising during its execution.\\

\subsection{Conjugate gradient method as a collection\\ of matrix graph traversals and error correlation}

We have shown above that, in the case $\mathrm{supp}(b)\!=\!\{i\}$, CG iterations are equivalent to a traversal of the matrix graph starting at vertex $i$. For a right-hand side $b$ with support of arbitrary cardinality, CG can be considered as a collection of traversals starting simultaneously from different vertices of the matrix graph $A$. These are the vertices $i\!\in\!V$ such that $i\!\in\!\mathrm{supp}(b)$. We show that the vector $d^k$ and each of its components $d_j^k$ can be decomposed into contributions that arrive at vertex $j$ from the vertices belonging to $\mathrm{supp}(b)$.

\smallskip

The following lemma holds.

\bigskip

\noindent \textbf{Lemma 1.} \textit{If $\mathrm{supp}(b)\!=\!\{i\}$ and $\ell(i,j)\!>\!k$, then $d_j^k\!=\!0$ and $\varepsilon_j^k\!=\!0$.}

\bigskip

\noindent \textbf{Proof.} We proceed by induction. Base case. For $k\!=\!1$, with $d^0\!=\!b$, the CG recurrence relations take the form 
$$d^1=r^1+\beta_1 d^0,$$ $$r^1=b-\alpha_0 Ab,$$ which gives $$d_j^1=(1+\beta_1)b_j-\alpha_0(Ab)_j.$$ If $\ell(i,j)\!>\!1$, then $a_{ji}\!=\!0$ and $b_j\!=\!0$, hence $$d_j^1\!=-\alpha_0\!(Ab)_j\!=-\alpha_0\!\sum\limits_{(j,t)\in E}^n a_{jt}b_t=-\alpha_0\cdot a_{ji}\cdot b_i=0.$$ 

For $k\!=\!2$: $$d_j^2=(1+\beta_2)d_j^1-\beta_1d_j^0-\alpha_1(Ad^1)_j.$$ If $\ell(i,j)\!>\!2$, then, assuming the statement holds for $k\!=\!1$, we have $d_j^1\!=\!0$ and $d_j^0\!=\!b_j\!=\!0$. Since $$(Ad^1)_j=\sum\limits_{(j,t)\in E}a_{jt}d_t^1\eqno(10)$$ and in the sum on the right-hand side of (10) the nonzero terms correspond only to those $t$ for which $d_t^1\!\neq\!0$. However, if $\ell(i,j)\!>\!2$, then for any such $t$ the corresponding factors~$a_{jt}$ in~(10) are equal to zero, since otherwise there would exist a walk $w\!=\!\{(i,t),(t,j)\}$ of length~$2$ in the matrix graph, which by the lemma's condition does not exist for $k\!=\!2$. Consequently, $d_j^2\!=\!0$. Thus, the statement of the lemma is proved for the cases $k\!=\!1$ and $k\!=\!2$.

Suppose that the statement of Lemma 1 holds for all values less than $k$. We show that it also holds for $k$. By the induction hypothesis, in (9) we have: $d_j^{k-1}\!=\!0$, $d_j^{k-2}\!=\!0$, and $$(Ad^{k-1})_j=\sum\limits_{(j,t)\in E}a_{jt}d_t^{k-1},$$ where in the sum on the right-hand side nonzero terms are possible only for those $t$ such that $d_t^{k-1}\!\neq\!0$. For such $t$, by the induction assumption, it must hold that $\ell(i,t)\!\le\!k\!-\!1$. If for such $t$ we have $a_{jt}\!\neq\!0$, this means that there exists a walk from $i$ to $j$ of length at most $k$. But by the lemma's condition $\ell(i,j)\!>\!k$. Hence, no such $a_{jt}$ exist, and $d_j^k\!=\!0$.

Since all components $d_j^{k-1}$, $d_j^{k-2}$, and $(Ad^{k-1})_j$ entering the expression (6) for $d_j^k$ are equal to zero, the corresponding operations in this expression also yield zero, and at the $k$-th iteration the error does not propagate to vertex $j$, nor does a round-off 
error arise in it. Consequently, $\varepsilon_j^k\!=\!0$. $\qed$

\bigskip

We introduce the following notation. Let
\begin{itemize}
\item[---] $d^k|_i$ denote the vector whose $j$-th component $d_j^k|_i$ is obtained by propagating the value of the $i$-th component $d_i^1$ of the vector $d^1$ over $k$ iterations;

\smallskip

\item[---] $\varepsilon^k|_i$ denotes the vector of computation errors of the vector $d^k$; each component $\varepsilon_j^k|_i$ of the vector $\varepsilon^k|_i$ is the error in computing the component $d_j^k$, which is obtained by propagating over $k$ iterations the value of $\varepsilon_i^1$---the $i$-th component of the error vector $\varepsilon^1$.
\end{itemize}

\bigskip

\noindent \textbf{Lemma 2.} $$d^k\!=\!\sum\limits_{i\in\mathrm{supp}(b)}d^k|_i,$$ $$\varepsilon^k\!=\!\sum\limits_{i\in\mathrm{supp}(b)}\varepsilon^k|_i.$$

\bigskip

\noindent \textbf{Proof.} To simplify the notation, consider the case $\mathrm{supp}(b)\!=\!\{1,2\}$: $$b\!=\!(b_1,b_2,0,\dots,0)\!=\! b|_1+b|_2,$$ 
where $b|_1\!=\!b_1\cdot e_1$, $b|_2\!=\!b_2\cdot e_2$. We prove the statement of the lemma for such $b$ by induction.

Base case. Let $d^0|_1\!=\! b|_1$ and $d^0|_2\!=\! b|_2$. For $k\!=\!1$, expression~(4) takes the form: $$d^1=(d^0|_1+d^0|_2)-\alpha_0A(d^0|_1+d^0|_2)=(d^0|_1-\alpha_0Ad^0|_1)+(d^0|_2-\alpha_0Ad^0|_2)=$$ $$=d^1|_1+d^1|_2.$$

For $k\!=\!2$: $$d^2=(1+\beta_2)d^1-\beta_1d^0-\alpha_1Ad^1=$$ $$=(1+\beta_2)(d^1|_1+d^1|_2)-\beta_1(d^0|_1+d^0|_2)-\alpha_1A(d^1|_1+d^1|_2)=$$
$$=\biggl((1+\beta_2)d^1|_1-\beta_1d^0|_1-\alpha_1Ad^1|_1\biggr)+\biggl ((1+\beta_2)d^1|_2-\beta_1d^0|_2-\alpha_1Ad^1|_2\biggr)=$$ $$=d^2|_1+d^2|_2.$$

Suppose that the statement of the lemma holds for iterations with numbers less than~$k$. We show that it also holds for the iteration with number $k$.

By the induction hypothesis, for iterations with numbers $k\!-\!2$ and $k\!-\!1$ we have: $$d^{k-2}=d^{k-2}|_1+d^{k-2}|_2,\ d^{k-1}=d^{k-1}|_1+d^{k-1}|_2,$$ hence, $$d^k=(1+\beta_k)d^{k-1}-\beta_{k-1}d^{k-2}-\alpha_{k-1}Ad^{k-1}=$$ $$=(1+\beta_k)(d^{k-1}|_1+d^{k-1}|_2)-\beta_{k-1}(d^{k-2}|_1+d^{k-2}|_2)-$$ $$-\alpha_{k-1}A(d^{k-1}|_1+d^{k-1}|_2)=\biggl((1+\beta_k)d^{k-1}|_1-\beta_{k-1}d^{k-2}|_1-\alpha_{k-1}Ad^{k-1}|_1\biggr)+$$ $$+\biggl ((1+\beta_k)d^{k-1}|_2-\beta_{k-1}d^{k-2}|_2-\alpha_{k-1}Ad^{k-1}|_2\biggr)=d^k|_1+d^k|_2.$$

We have shown that $d^k\!=\!d^k|_1+d^k|_2$ for $b\!=\!b|_1+b|_2$. Since for an arbitrary $b$ we have $$b\!=(b_1,\dots,b_n)=\!\sum\limits_{i\in\mathrm{supp}(b)}b|_i,$$ where $b|_i\!=\!b_i\cdot e_i$, and since the recurrence relation (6) is linear, we obtain $$d^k=\sum\limits_{i\in\mathrm{supp}(b)}d^k|_i.$$

Since for $\varepsilon^k$ expression~(9) is completely analogous to~(6), the equality $$\varepsilon^k=\sum\limits_{i\in\mathrm{supp}(b)}\varepsilon^k|_i$$ is proved similarly. $\qed$

\bigskip

From Lemmas 1 and 2 it follows

\bigskip

\noindent\textbf{Lemma 3.} $$d_j^k=\sum\limits_{\{i\, |\, \ell(i,j)\le k\}}d_j^k|_i,$$ $$\varepsilon_j^k=\sum\limits_{\{i\, |\,\ell(i,j)\le k\}}\varepsilon_j^k|_i.$$

\bigskip

\noindent \textbf{Proof.} By Lemma~1, if $\ell(i,j)\!>\!k$, then for the right-hand side $b|_i\!=\!b_i\cdot e_i$, we have $$d_j^k=d_j^k|_i=0,\quad \varepsilon_j^k=\varepsilon_j^k|_i=0.$$ Since by Lemma 2 $$d_j^k=\sum\limits_{i\in\mathrm{supp}(b)}d_j^k|_i,\quad\varepsilon_j^k=\sum\limits_{i\in\mathrm{supp}(b)}\varepsilon_j^k|_i,$$ and $\mathrm{supp}(b)=\{i\, |\, \ell(i,j)\!\le\!k\}\bigcup\{i\, |\, \ell(i,j)\!>\!k\}$, the statement of Lemma~3 is proved. $\qed$

\bigskip

Thus, the information obtained at the first iteration, namely the components of the vector $d^1$, including the error $\varepsilon^{1}$ arising from round-off errors at the first iteration, propagates only between those vertices of the matrix graph whose distance does not exceed $k$. Consequently, 
errors in two different vertices of the graph (components of the search direction vectors) do not affect each other, that is, they are uncorrelated, if the number of CG iterations performed is less than the distance between them. This follows from Theorem~2.

\bigskip

\noindent \textbf{Theorem 2.} \textit{Errors caused by round-off errors in the $i$-th component of the search direction vector at the $k_1$-th CG iteration affect the error in computing the $j$-th component of the search direction vector, that is, give a nonzero contribution to it, 
not earlier than at the iteration with number $k_2\!=\!k_1+\ell(i,j)$.}

\bigskip

Theorem 2 follows from Lemmas~1--3 for the case $k_1\!=\!1$. For an arbitrary $k_1$ this holds, since the round-off error arising at the $k_1$-th iteration in the $i$-th component propagates along the matrix graph in exactly the same way as the error arising at the first iteration, but with a shift of $k_1\!-\!1$ iterations. Therefore, it will reach vertex $j$ no earlier than at iteration $k_2\!=\!k_1+\ell(i,j)$.

\bigskip

Since $$x^k=x^{k-1}+\alpha_{k-1}d^{k-1},\quad r^k=r^{k-1}-\alpha_{k-1}Ad^{k-1},$$ the errors in computing the components of the approximate solution and the residual also have time to propagate by the end of the $k$-th iteration between the vertices of the matrix graph that are at a distance not exceeding $k$ from each other: $$x_j^k=\sum\limits_{\{i\, |\, \ell(i,j)\le k\}}x_j^k|_i,\quad r_j^k=\sum\limits_{\{i\, |\,\ell(i,j)\le k\}}r_j^k|_i,$$
where $x_j^k|_i$ and $r_j^k|_i$ are the values of the components of the vectors $x^k$ and $r^k$ obtained by propagating over $k$ CG iterations the value of the $i$-th component of the vector $d^1$. Similar decompositions hold for the errors of these vectors as well.\\

\section{Walk contributions to the error}

The decomposition of the error obtained at the $k$-th iteration into contributions of walks of length not exceeding $k$ follows from the linearity of the CG recurrence relations (5), (8) and our additive model of round-off errors. From (5) and (8) it follows that contributions are delivered to both 
the components of the search direction vector and the corresponding components of the error vector along the same walks in the matrix graph. By Theorem~2, if $\ell(i,j)\!=\!k$, then at the $k$-th iteration we have $d_j^{k-1}\!=\!d_j^{k-2}\!=\!0$ and $\varepsilon_j^{k-1}\!=\!\varepsilon_j^{k-2}\!=\!0$. 
Consequently, when at the $k$-th iteration the contribution from vertex $i$ reaches vertex $j$, we have: $$d_j^k=-\alpha_{k-1}\sum\limits_{(j,t)\in E}a_{jt}d_t^{k-1},\eqno(11)$$ $$\varepsilon_j^k=-\alpha_{k-1}\sum\limits_{(j,t)\in E}a_{jt}\varepsilon_t^{k-1}.\eqno(12)$$ From this it follows that after performing the $k$-th iteration, similarly to the decomposition into contributions of information transmitted to vertex $j$ during the BFS traversal (3), scaled by the factor $-\alpha_{k-1}$, for a vector $b$ with arbitrary support we obtain: $$d_j^k=\sum\limits_{i\in\mathrm{supp}(b)}\sum\limits_{w\in 
W_{ij}^k}C(w),\eqno(13)$$ where for a walk $w\!=\!\{(i_0,i_1),\dots,(i_{k-1},i_k)\}\!\in\!W_{ij}^k$ with $i_0\!=\!i$, $i_k\!=\!j$, we obtain that its contribution $C(w)$ is determined by the value $d_i^0$ and the \textit{weight of the walk}: $$C(w)=\biggl (\prod_{q=1}^k(-\alpha_{q-1}a_{i_{q-1}i_q})\biggr )d_i^0.$$ The product $$v(w)\!=\!\prod_{q=1}^k(-\alpha_{q-1}a_{i_{q-1}i_q})\eqno(14)$$ is the weight of the walk $w$. At the $(k\!+\!1)$-th 
iteration, the information obtained during the algorithm initialization will propagate further to vertices adjacent to~$j$.

Similarly to (13), the error $\varepsilon_i^1$ arising at the first iteration propagates along a walk $w$ of length $k$ starting at vertex $i$ and ending at vertex $j$. Such a walk contributes $C_e(w)$ to the error $\varepsilon_j^k$ in computing $d_j^k$ at the $k$-th iteration: $$C_e(w)=v(w)\,\varepsilon_i^1.$$

After the $k$-th iteration, the nonzero terms $(1+\beta_t)\varepsilon_j^{t-1}$ and $-\beta_{t-1}\varepsilon_j^{t-2}$ for $t\!=\!k+1, k+2$ will contribute nonzero values to the errors $\varepsilon_j^{k+1}$ and $\varepsilon_j^{k+2}$. Note that in exact arithmetic, as in our additive model of 
information propagation during CG, the coefficients $\beta_k$ monotonically tend to zero as $k$ increases, and consequently, the processes of information propagation~(6) and error propagation~(9) at later iterations increasingly closely correspond to~(11) and~(12).

Thus, the error $\varepsilon_i^1$ arising at the first iteration propagates along walks in the graph. The magnitude of its contribution $C_e(w)$, delivered along walk $w$ to the total error in computing the component of the search direction vector at the $k$-th iteration, is determined by its 
weight~(14). Errors caused by round-off errors at subsequent iterations propagate along the same walks with the same weight.

The order of magnitude of the error $\varepsilon_i^1$ in computing the $i$-th component of the search direction vector at the first iteration is $O(\varepsilon_{\mathrm{mach}})$. When this error propagates to vertex $j$, that is, to the $j$-th component of the search direction vector at the $k$-th 
iteration, its value is multiplied by the weight $v(w)$. In connected sparse graphs, even in a path graph on $n$ vertices, the number of walks along which errors propagate grows exponentially with $k$, and hence the number of contributions to each component of the error vector grows exponentially. 
Although this growth may be balanced by the exponential decay of walk weights or by the alternation of signs of contributions, the correlation of errors leads to stagnation of the iterative process.\\

\subsection{Loss of orthogonality of residuals and conjugacy of directions due to accumulation and correlation of round-off errors}

The residuals $\hat{r}^{k_1}$, $\hat{r}^{k_2}$ and the search direction vectors $\hat{d}^{k_1}$, $\hat{d}^{k_2}$ computed in exact arithmetic at the $k_1$-th and $k_2$-th iterations satisfy the following conditions: $$(\hat{r}^{k_1},\hat{r}^{k_2})=0\mbox{\ and\ }(\hat{d}^{k_1},A\hat{d}^{k_2})=0.$$ 
For computations in floating-point arithmetic $$(r^{k_1},r^{k_2})=(\hat{r}^{k_1}+\dot{\epsilon}^{k_1},\hat{r}^{k_2}+\dot{\epsilon}^{k_2})=(\hat{r}^{k_1},\hat{r}^{k_2})+(\hat{r}^{k_1},\dot{\epsilon}^{k_2})+(\dot{\epsilon}^{k_1},\hat{r}^{k_2})+(\dot{\epsilon}^{k_1},\dot{\epsilon}^{k_2})=$$ 
$$=(\hat{r}^{k_1},\dot{\epsilon}^{k_2})+(\dot{\epsilon}^{k_1},\hat{r}^{k_2})+(\dot{\epsilon}^{k_1},\dot{\epsilon}^{k_2}),\eqno(15)$$ where $\dot{\epsilon}^{k_1}$ and $\dot{\epsilon}^{k_2}$ are the errors accumulated by the $k_1$-th and $k_2$-th iterations in computing $r^{k_1}$ and~$r^{k_2}$. Similarly: $$(d^{k_1},Ad^{k_2})=\bigl((\hat{d}^{k_1}+\epsilon^{k_1}),A(\hat{d}^{k_2}+\epsilon^{k_2})\bigr)=$$ $$=(\hat{d}^{k_1},A\hat{d}^{k_2})+(\hat{d}^{k_1},A\epsilon^{k_2})+(\epsilon^{k_1},A\hat{d}^{k_2})+(\epsilon^{k_1},A\epsilon^{k_2})=$$ $$=(\hat{d}^{k_1},A\epsilon^{k_2})+(\epsilon^{k_1},A\hat{d}^{k_2})+(\epsilon^{k_1},A\epsilon^{k_2}),\eqno(16)$$ where $\epsilon^{k_1}$ and $\epsilon^{k_2}$ are, respectively, the vectors of accumulated errors in computing $d^{k_1}$ and $d^{k_2}$ by the $k_1$-th and $k_2$-th iterations.

The paths traversed by the errors arising at early iterations are traversed also by the errors arising at later ones, giving on average approximately the same contributions. The emerging dependencies between the error components at different iterations create a systematic correlation between the errors at different iterations. For this reason, the loss of orthogonality and conjugacy occurs significantly earlier than in~the case of independent errors $\epsilon^{k_1}$ and $\epsilon^{k_2}$. The first-order terms in (15) and (16) make the main contribution to the loss of orthogonality, evaluated through the non-equality to zero of the corresponding scalar products, but this loss is conditioned not only by the absolute magnitude of the errors, but also by their mutual correlation. 

In the case of independence of errors at the $k_1$-th and $k_2$-th iterations, the sum of these terms is on average equal to zero, just as the second-order terms $(\dot{\epsilon}^{k_1},\dot{\epsilon}^{k_2})$ and $(\epsilon^{k_1},A\epsilon^{k_2})$ are on average equal to zero. Starting from~some iteration, the correlation of errors becomes increasingly global, such compensation does not occur, and~the scalar products $(r^{k_1},r^{k_2})$ and~$(d^{k_1},Ad^{k_2})$ become increasingly distinct from zero.\\

\section{Diameter $\mathcal{D}$, effective diameter $\mathcal{D}_{\mathrm{eff}}$\\ and pseudo-diameter $\tilde{\mathcal{D}}$ of the matrix graph}

In the experiments below, we record the discrepancy of the residuals for computations in single and double precision, and consider the dependence of the iteration $N_{\mathrm{div}}$, at which it was recorded, on the diameter, effective diameter, and pseudo-diameter of the matrix graph. The effective diameter $\mathcal{D}_{\mathrm{eff}}$ is a more accurate indicator of the accumulation and amplification of error correlation than the diameter $\mathcal{D}$, since for matrices in which the majority of vertices are connected by relatively short paths, when calculating $\mathcal{D}_{\mathrm{eff}}$, rare long paths are not taken into account. Such paths can be considered as outliers when estimating the average distance between vertices.

For matrices with extended star graphs and Laplacian matrices of differential equations considered below, both the values of the diameter $\mathcal{D}$ and the effective diameter $\mathcal{D}_{\mathrm{eff}}$ correlate equally well with the iteration number $N_{\mathrm{div}}$, at which the divergence of residuals occurs. This is explained by the fact that for almost all vertices in these graphs, the most distant vertices are located at comparable distances, and the difference in the values of $\mathcal{D}$ and $\mathcal{D}_{\mathrm{eff}}$ is not critical for estimating the relationship between the values of $N_{\mathrm{div}}$ for matrices of these types.

For random matrices, whose graphs are generated according to the Erdős-Rényi model \cite{Erdos}, where each pair of vertices is connected by an edge independently of other pairs with probability $p$, the difference between $\mathcal{D}$ and $\mathcal{D}_{\mathrm{eff}}$ can be large, since for the majority of vertices, the vertices most distant from them are located at a distance of order $O(\log n)$. Therefore, the dependence between $\mathcal{D}_{\mathrm{eff}}$ and $N_{\mathrm{div}}$ is more pronounced than between $\mathcal{D}$ and $N_{\mathrm{div}}$. However, in the experiments with random matrices, we will consider the dependence between $\tilde{\mathcal{D}}$ and $N_{\mathrm{div}}$. The pseudo-diameter $\tilde{\mathcal{D}}$ is used by us in these experiments for two reasons: first, the repeated computation of $\tilde{\mathcal{D}}$ for all graphs in the sample is performed significantly faster than the computation of $\mathcal{D}_{\mathrm{eff}}$; second, this shows that for two matrices of the same type, even the closeness of such a rough metric characteristic of the matrix graph, as the pseudo-diameter, corresponds to the closeness of the values of $N_{\mathrm{div}}$ for linear systems of algebraic equations with these matrices.

Therefore, in the experiments considered below for matrices with extended star graphs and Laplacian matrices of partial differential equations, we analyze the dependence between $\mathcal{D}$, $\mathcal{D}_{\mathrm{eff}}$, and $N_{\mathrm{div}}$, and for random sparse matrices---between $\tilde{\mathcal{D}}$ and $N_{\mathrm{div}}$.\\

\section{Computational Experiments}

The experiments considered below demonstrate the dependence of the iteration number $N_{\mathrm{div}}$, at which the divergence of residuals is recorded, on the metric characteristics of the matrix graph. For computations in single precision, we denote this number as $N_{\mathrm{div}}^{\mathrm{\textit{sp}}}$, and in double precision---as~$N_{\mathrm{div}}^{\mathrm{\textit{dp}}}$.

\subsection{Matrix with a star graph and matrix with a path graph}

\begin{figure}[htbp]
\begin{center}
\includegraphics[width=130mm]{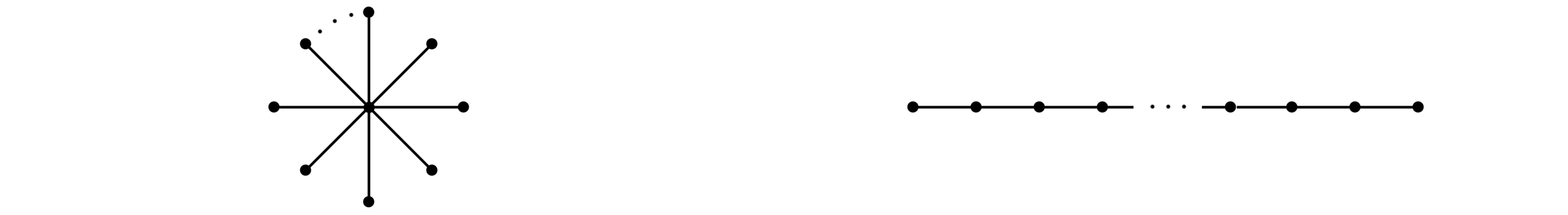}
\caption{Star and path graphs.}\label{fig:5}
\end{center}
\end{figure}

Consider matrices obtained by modifying the adjacency matrices of star and path graphs on $601$ vertices (Fig.~~\ref{fig:5}). The diagonal elements of these matrices are chosen so as to ensure strict diagonal dominance:
$$a_{ii}=\sum\limits_{j=1}^na_{ij}+\delta.\eqno(17)$$
Nonzero diagonal elements correspond to weighted loops at the vertices of the graph. The value $\delta\!>\!0$ in (17) allows one to control the condition number of the matrix. In what follows, we set $\delta\!=\!0.1$.

For the matrix with a star graph, the diameter $\mathcal{D}$ takes the smallest possible value for a connected graph on $n$ vertices: $\mathcal{D}\!=\!2$. Whereas for the matrix with a path graph, it takes the largest possible value: $\mathcal{D}\!=\!n\!-\!1$, which in the case under consideration gives $\mathcal{D}\!=\!600$.

\begin{figure}[htbp]
\begin{center}
\includegraphics[width=110mm]{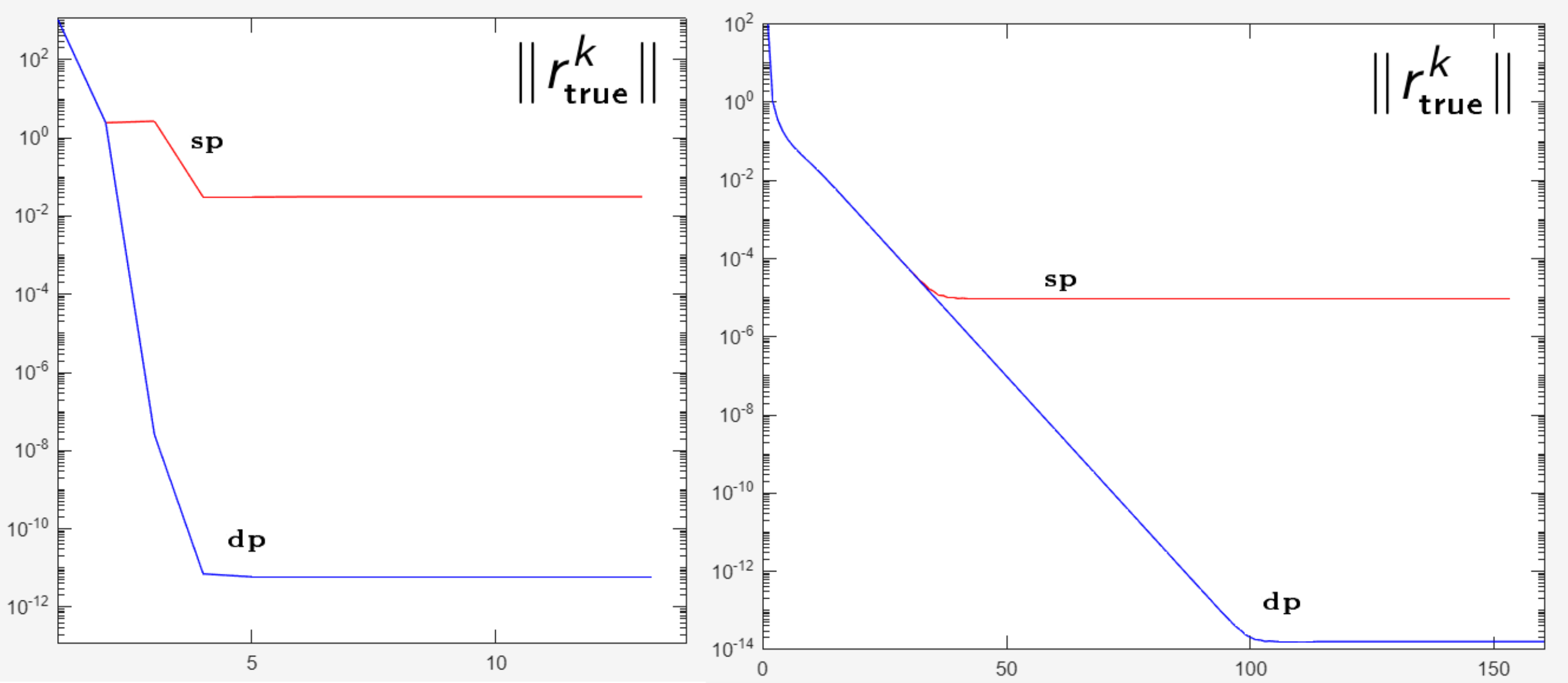}
\caption{True residuals for matrices with a star graph (left, $\mathcal{D}\!=\!2$) and a path graph (right, $\mathcal{D}\!=\!600$) in two precision modes.}\label{fig:6}
\end{center}
\end{figure}

As shown in Fig.~\ref{fig:6}, for the matrix with a star graph, stagnation begins at the 4th iteration in computations in both double and single precision. For the matrix with a path graph, $35$ iterations are performed before the onset of the stagnation stage in single precision, and $100$ iterations in double precision. That is, for the matrix with $\mathcal{D}\!=\!2$, the propagation of round-off errors leads to~stagnation an order of magnitude or more earlier than for the matrix with $\mathcal{D}\!=\!600$.\\

\subsection{Matrices with extended star graphs}

\begin{figure}[htbp]
\begin{center}
\includegraphics[width=40mm]{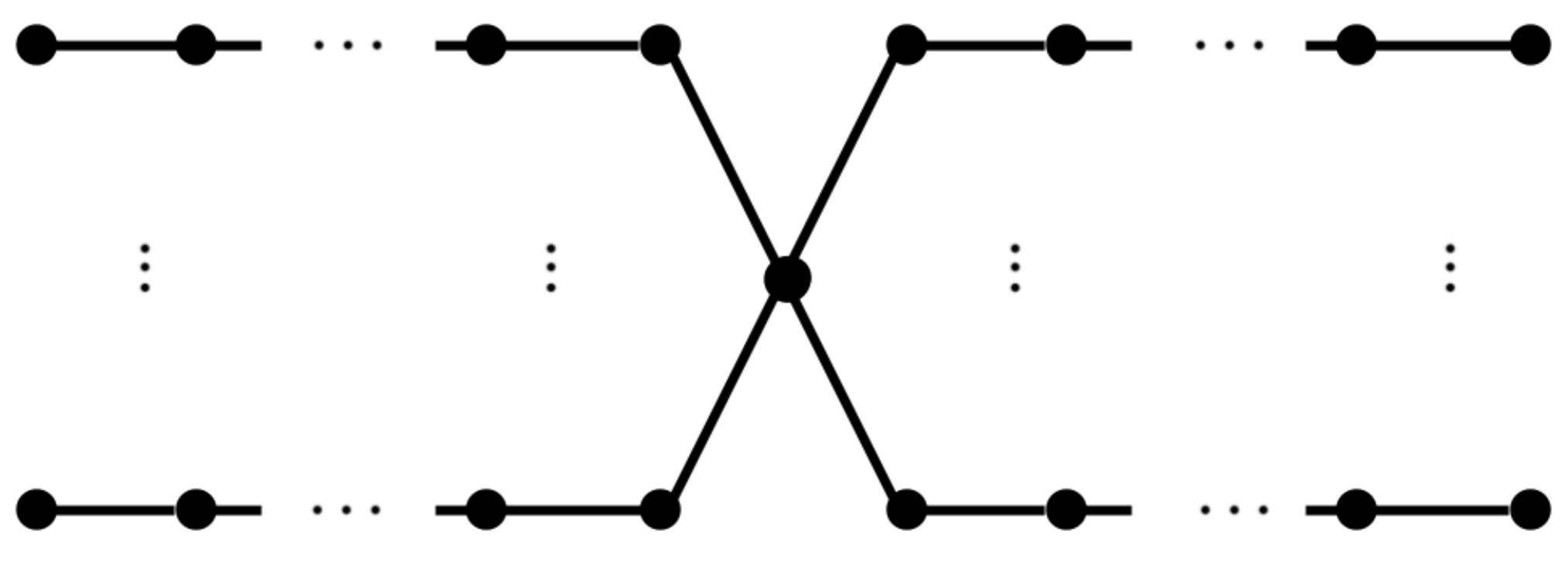}
\caption{Extended star graph.}\label{fig:7}
\end{center}
\end{figure}

Consider matrices obtained by introducing diagonal dominance into the adjacency matrices of \textit{extended stars} in accordance with~(17)---graphs in which one central vertex has paths of the same length $l$ attached to it (Fig.~\ref{fig:7}). An extended star is defined by two parameters $(k_{\mathrm{p}}, l_{\mathrm{p}})$, where $k_{\mathrm{p}}$ is the number of paths connected to the central vertex, and $l_{\mathrm{p}}$ is their length.

Table~1 presents the condition number $\kappa$, the diameter $\mathcal{D}$, the effective diameter $\mathcal{D}_{\mathrm{eff}}$, the iteration number $N_{\mathrm{div}}^{\mathrm{\textit{\textit{sp}}}}$ at which the divergence of residuals occurs in single-precision computations, and the iteration number $N_{\mathrm{div}}^{\mathrm{\textit{dp}}}$ at which it occurs in double-precision computations. As can be seen, in both computation modes, the divergence of residuals occurs the later, the larger the diameter of the matrix graph.\\

\begin{table}[H]
\begin{center}
\caption{Iteration number $N_{\mathrm{div}}$, at which the divergence of residuals occurs for~equally conditioned matrices with graphs of $(k_{\mathrm{p}}, l_{\mathrm{p}})$-extended stars. $n\!=\!601$.}

\bigskip

\begin{tabular}{|c|c|c|c|c|c|c|c|c|c|}
\hline  
\quad             & $(600,1)$       & $(120,5)$         & $(20,30)$        & $(10,60)$        & $(3,200)$    \\
\hline
\hline
$\kappa$          & $6011.0$  & $1211.1$  & $311.4$  & $112.1$  & $46.0$\\
\hline 
\hline
$\mathcal{D}$               & $2$      & $10$     & $40$     & $120$    & $400$\\
\hline
$\mathcal{D}_{\mathrm{eff}}$             & $2$      & $9$      & $32$     & $93$     & $291$\\
\hline
\hline
$N_{\mathrm{div}}^{\mathrm{\textit{sp}}}$          & $4$      & $12$      & $28$     & $41$     & $45$\\
\hline
$N_{\mathrm{div}}^{\mathrm{\textit{dp}}}$          & $4$      & $9$      & $25$     & $66$     & $105$\\
\hline
\end{tabular}  
\end{center}  
\end{table}

\begin{table}[H]
\begin{center}
\caption{Iteration number $N_{\mathrm{div}}$, at which the divergence of residuals occurs for~equally conditioned matrices with graphs of $(k_{\mathrm{p}}, l_{\mathrm{p}})$-extended stars. $n\!=\!8001$.}

\bigskip

\begin{tabular}{|c|c|c|c|c|c|c|c|c|}
\hline  
\quad  & (8000,1) & (2000,4) & (800,10) & (200,40) & (40,200) & (10,800) & (5,1600) & (2,4000) \\
\hline
\hline

$\mathcal{D}$ & $2$      &  $8$     &  $20$    &  $80$    &  $400$   & $1600$   &  $3200$ &  $8000$\\
\hline
$\mathcal{D}_{\mathrm{eff}}$ &   $2$    &   $7$   &   $16$    &   $63$    &  $310$    &  $1224$   &  $2401$  &  $5471$ \\
\hline
\hline
$\kappa$  & \multicolumn{8}{|c|}{$N_{\mathrm{div}}^{\mathrm{\textit{sp}}}$} \\
\hline
$10^2$  & $4$      &  $5$     &  $7$    &  $11$    &  $21$    & $37$     &  $44$ &  $55$\\
\hline
$10^3$  & $4$      &  $5$     &  $10$    &  $27$    &  $57$    & $94$     &  $106$ &  $131$\\
\hline
$10^4$  & $4$      &  $6$     &  $12$    &  $36$    &  $131$    & $185$     &  $216$ &  $307$\\
\hline
$10^5$  & $4$      &  $9$     &  $16$    &  $45$    &  $160$    & $451$     &  $495$ &  $564$\\
\hline
\hline
$\kappa$  & \multicolumn{8}{|c|}{$N_{\mathrm{div}}^{\mathrm{\textit{dp}}}$} \\
\hline
$10^2$  & $4$      &  $7$     &  $12$    &  $24$    &  $45$    & $83$   &  $107$ &  $152$\\
\hline
$10^3$  & $4$      &  $7$     &  $14$    &  $48$    &  $116$    & $222$   &  $305$ &  $443$\\
\hline
$10^4$  & $4$      &  $8$     &  $16$    &  $50$    &  $226$    & $612$     &  $847$ &  $1219$\\
\hline
$10^5$  & $4$      &  $9$     &  $17$    &  $51$    &  $227$    & $855$     &  $1658$ &  $3407$\\
\hline
\end{tabular}  
\end{center}  
\end{table}

\subsection{Equally conditioned matrices with extended star graphs}

By choosing $\delta$ in~(17) so that the matrix has a given condition number, identical for all matrices, we obtain a more clearly expressed dependence between the diameter of the matrix graph and the iteration at which the divergence of residuals is recorded. Table~2 shows the values of $N_{\mathrm{div}}^{\mathrm{\textit{sp}}}$ and $N_{\mathrm{div}}^{\mathrm{\textit{dp}}}$ for matrices with extended star graphs on~$8001$ vertices with condition numbers equal to $10^2$, $10^3$, $10^4$, and $10^5$.

The data in Table~2 show that with an increase in the condition number, the divergence of residuals for matrices with larger diameter is recorded later than for~matrices with~the same condition number, but a smaller diameter. That is, a larger condition number enhances the sensitivity of the value $N_{\mathrm{div}}$ to differences in $\mathcal{D}$ (or $\mathcal{D}_{\mathrm{eff}}$): the~greater the difference in the diameters of the matrix graphs, the greater the differences in the values of $N_{\mathrm{div}}$ for~them.\\

\subsection{Laplacian matrices of partial differential equations}

The Laplace operator $\Delta\!=\!\frac{\partial^2}{\partial x^2}+\frac{\partial^2}{\partial y^2}+\frac{\partial^2}{\partial z^2}$ is one of the fundamental objects of mathematical physics. It appears in key partial differential equations: the Poisson equation $\Delta u\!=\!f$ (electrostatics, stationary heat conduction, theory of elasticity), the heat equation $\frac{\partial u}{\partial t}\!=\!\alpha\Delta u$, and the Helmholtz equation $\Delta u+k^2u\!=\!f$. When solving these equations numerically using the finite difference method, the continuous domain is replaced by a discrete grid of $n\!=\!n_1\!\times\!n_2\!\times\!n_3$ nodes, and the action of the Laplace operator on a grid function is approximated by a matrix $L$ of size $n\!\times\!n$.

The discrete Laplacian matrix is highly sparse: in the three-dimensional case, there are no more than $7$ non-zero elements per node (the value at the center of the stencil and $6$ neighbors along the coordinate axes), therefore, the proportion of non-zero elements is of order $O(1/n)$ and~tends to zero as the grid is refined. We used Laplacian matrices arising from problems with homogeneous Dirichlet boundary conditions, under which $L$ is positive definite.

The condition number $\kappa(L)$ depends on both the grid size and its geometry. In~particular, for grids of a fixed volume $n\!=\!n_1\!\times\!n_2\!\times\!n_3$, elongated configurations (with a large effective diameter $\mathcal{D}_{\mathrm{eff}}\!\approx\!\sqrt{n_1^2+n_2^2+n_3^2}$) 
lead to significantly larger values of $\kappa$ than compact, cube-like configurations.

\begin{figure}[H]
\centering
\includegraphics[width=120mm]{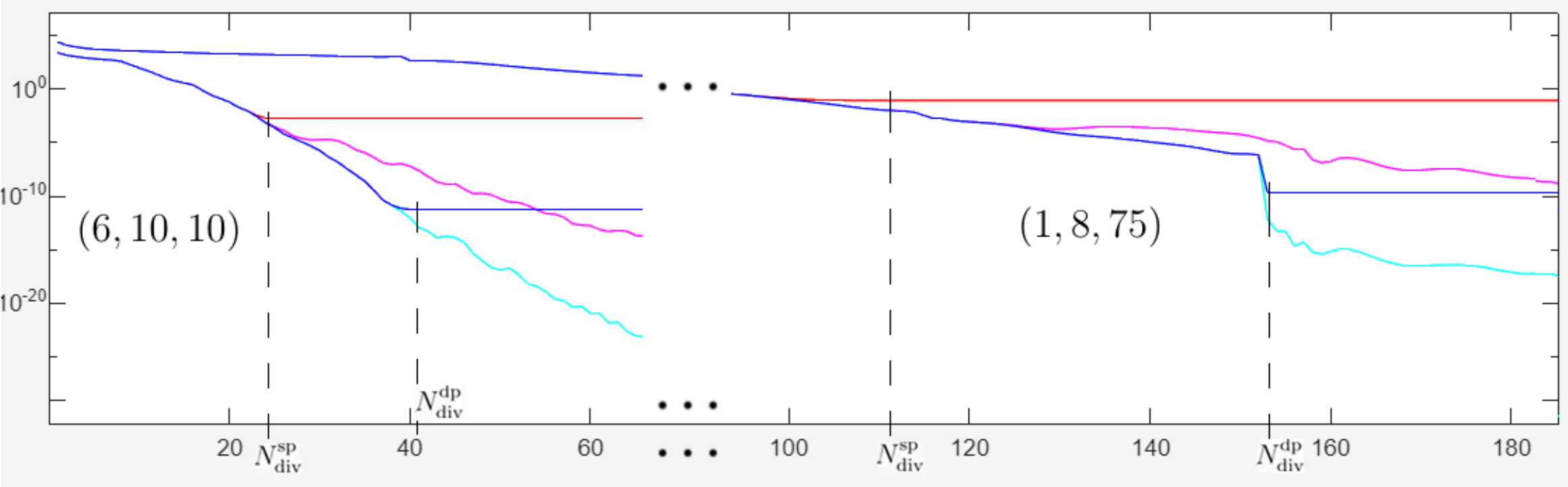}
\caption{Dependence of residuals on the CG iteration number in double-precision computations for Laplacians with different diameters (grid configurations---$(6,10,10)$ and $(1,8,75)$).}\label{fig:8}
\end{figure}

Fig.~\ref{fig:8} presents the graphs of the true and recursive residuals for Laplacians with configurations $(n_1,n_2,n_3)\!=\!(6,10,10)$ and $(n_1,n_2,n_3)\!=\!(1,8,75)$, having diameters $23$ and $81$, respectively. These graphs show that, although in the case of a larger diameter the convergence rate is lower, the divergence of residuals occurs later. The values in Table~3 show that the dependence between the iteration $N_{\mathrm{div}}$, at which the divergence of residuals is recorded, and the diameter of the matrix graph, observed for extended stars, also holds for Laplacians with different diameters.

\bigskip

\noindent\textbf{Remark 2.} Despite the fact that for matrices with extended star graphs the condition number decreases with the diameter (Table~1), whereas for Laplacians, on the contrary, it increases (Table~3), in both cases $N_{\mathrm{div}}$ increases with the diameter. In other words, when these two cases are considered jointly, no correlation is observed between such a spectral characteristic as the condition number and $N_{\mathrm{div}}$. At the same time, there exists a clear correlation between the metric characteristics of the matrix graph---the diameter $\mathcal{D}$ and the effective diameter $\mathcal{D}_{\mathrm{eff}}$---and the divergence of residuals; moreover, the dependence of $N_{\mathrm{div}}$ on $\mathcal{D}$ and $\mathcal{D}_{\mathrm{eff}}$ is direct (positive).

\begin{table}[H]
\begin{center}
\caption{The iteration number $N_{\mathrm{div}}$, at which the divergence of residuals occurs, for~Laplacians in two precision modes of computations.}

\bigskip

\begin{tabular}{|c|c|c|c|c|c|c|c|c|c|}
\hline  
$(n_1,n_2,n_3)$ & $(6,10,10)$ & $(3,8,25)$ & $(1,15,40)$& $(1,8,75)$\\
\hline
$\kappa$        & $62.49$  & $166.02$   &  $392.27$ &  $1170.51$\\
\hline 
$\mathcal{D}$             &  $23$     &  $33$      & $53$     & $81$\\
\hline
$\mathcal{D}_{\mathrm{eff}}$           &  $13$     &  $21$     & $33$     & $54$\\
\hline
\hline
$N_{\mathrm{div}}^{\mathrm{\textit{sp}}}$          & $26$      & $58$      & $65$     & $113$\\
\hline
$N_{\mathrm{div}}^{\mathrm{\textit{dp}}}$          & $41$      & $94$      & $102$     & $153$\\
\hline
\end{tabular}  
\end{center}  
\end{table}

\subsection{Random matrices}

Consider matrices generated as follows.

\begin{itemize}
\item[1.] A sparse random graph $G$ with an adjacency matrix $A^0\!=\!(a_{ij}^0)$ is generated.

\item[2.] Based on the obtained matrix $A^0$, the matrix $A$ is formed:

\begin{itemize}
\item[2.1] For each $a_{ij}^0\!\neq\!0$, a value $a_{ij}\!\in\![-10,10]$ is generated according to a uniform distribution on the interval $[-10,10]$.

\smallskip

\item[2.2] The diagonal elements $a_{ii}$ are formed in accordance with~(17), choosing the value of $\delta$ for the matrix $A$ such that $A$ has a prescribed condition number.
\end{itemize}
\end{itemize}

In the computational experiment, the results of which are presented in Table~4, for each value of $\kappa\!=\!10^3, 10^4, 10^5$, $100$ equally conditioned random matrices of dimension $n\!=\!1000$ were generated, with mean values of the diameter close to~$100$ and $900$ (the standard deviation of the pseudo-diameter in each of these cases did not exceed $6\%$ of these values of $\mathcal{D}$). The mean values $\mu(N_{\mathrm{div}})$ and the corresponding standard deviations $\sigma(N_{\mathrm{div}})$ are presented for~computations in~two precision modes. Table~5 presents the results of the same experiment for similarly generated matrices of dimension $n\!=\!10000$.

\begin{table}[H]
\begin{center}
\caption{Mean value and standard deviation of $N_{\mathrm{div}}$ for random matrices with~given diameters in two precision modes. $n\!=\!1000$.}

\bigskip

\begin{tabular}{|c||c|c||c|c||c|c|}
\hline
$\kappa$ & \multicolumn{2}{|c||}{$10^3$} & \multicolumn{2}{|c||}{$10^4$} & \multicolumn{2}{|c|}{$10^5$}\\
\hline
$\tilde{\mathcal{D}}$ & $\approx 100$ & $\approx 900$ & $\approx 100$ & $\approx 900$ & $\approx 100$ & $\approx 900$\\
\hline
\hline
$\mu\bigl(N_{\mathrm{div}}^{\mathrm{\textit{sp}}}\bigr)$ & 101.53 & 192.13 & 438.20 & 484.53 & 947.11 & 1086.17\\
\hline
$\mu\bigl(N_{\mathrm{div}}^{\mathrm{\textit{dp}}}\bigr)$ & 402.32 & 455.07 & 1112.11 & 1265.22 & 1871.38 & 2006.09\\
\hline
\hline
$\sigma\bigl(N_{\mathrm{div}}^{\mathrm{\textit{sp}}}\bigr)$  & 5.72   &  5.88 & 15.23   & 13.65 & 46.35   & 40.65\\
\hline
$\sigma\bigl(N_{\mathrm{div}}^{\mathrm{\textit{dp}}}\bigr)$ & 16.08 & 16.49 & 41.26 & 38.50 & 83.89 & 90.13\\
\hline
\end{tabular}

\end{center}
\end{table}

It is evident that $\mu(N_{\mathrm{div}}^{\mathrm{\textit{sp}}})$ and $\mu(N_{\mathrm{div}}^{\mathrm{\textit{dp}}})$ for the generated random matrices increase with an increase in the diameter of the matrix graph. The difference in the values of $N_{\mathrm{div}}$ between matrices with different pseudo-diameters turns out to be somewhat less pronounced compared to matrices with extended star graphs and Laplacian matrices. This difference persists across all pairs of prescribed values of $\mathcal{D}$, and in the majority of cases, it increases with the condition number.

\begin{table}[H]
\begin{center}
\caption{Mean value and standard deviation of $N_{\mathrm{div}}$ for random matrices with~given diameters in two precision modes. $n\!=\!10000$.}

\bigskip

\begin{tabular}{|c||c|c||c|c||c|c|}
\hline
$\kappa$ & \multicolumn{2}{|c||}{$10^3$} & \multicolumn{2}{|c||}{$10^4$} & \multicolumn{2}{|c|}{$10^5$}\\
\hline
$\tilde{\mathcal{D}}$ & $\approx 1000$ & $\approx 9000$ & $\approx 1000$ & $\approx 9000$ & $\approx 1000$ & $\approx 9000$\\
\hline
\hline
$\mu\bigl(N_{\mathrm{div}}^{\mathrm{\textit{sp}}}\bigr)$  & 183.00 & 199.40 & 468.30 & 494.33 & 1132.20 & 1212.10\\
\hline
$\mu\bigl(N_{\mathrm{div}}^{\mathrm{\textit{dp}}}\bigr)$ &  431.20 & 479.20 & 1252.10 & 1271.52 & 3609.50 & 3961.60\\
\hline
\hline
$\sigma\bigl(N_{\mathrm{div}}^{\mathrm{\textit{sp}}}\bigr)$  & 2.35   &  6.99 & 11.31   & 12.65 & 29.79   & 31.83\\
\hline
$\sigma\bigl(N_{\mathrm{div}}^{\mathrm{\textit{dp}}}\bigr)$  & 6.46    & 17.48 & 41.26 & 36.43  &  115.89 & 127.23\\
\hline
\end{tabular}

\end{center}
\end{table}

In addition, this experiment with random matrices shows that close values of the pseudo-diameters of two graphs correspond to close values of their $N_{\mathrm{div}}$. This is confirmed by the values of $\sigma\bigl(N_{\mathrm{div}}\bigr)$, which demonstrate that for random matrices with given diameters, the values of $N_{\mathrm{div}}$ exhibit low variability, as evidenced by the small standard deviation. In most cases presented in Tables~4 and 5, it remains below $5\%$ of the mean value of $N_{\mathrm{div}}$. With an increase in the matrix size, the spread of the values of $N_{\mathrm{div}}$, evaluated via $\sigma\bigl(N_{\mathrm{div}}\bigr)$, for matrices with close $\tilde{\mathcal{D}}$ becomes smaller relative to the mean value (Table~5).\\

\section{Using the metric characteristics of the matrix graph\\ to switch the computation precision mode}

Mixed-precision computation algorithms~\cite{Higham} are based on switching precision modes when solving problems of numerical linear algebra 
(solving systems of linear equations, matrix factorizations, etc.). As an example of a mixed-precision computation algorithm for solving systems of linear equations, consider the following two-stage computation scheme. The computations are performed in two precision modes: at the first stage---in single precision, at the second stage---in double precision. 

\begin{itemize}
\item[I.] Perform $N_{\mathrm{\textit{sp}}}$ CG iterations in single precision, obtaining an approximate solution~$x_{\mathrm{\textit{dp}}}^0$.
\item[II.] Using $x_{\mathrm{\textit{dp}}}^0$ as the initial approximation, perform $N_{\mathrm{\textit{dp}}}$ CG iterations in double precision, obtaining an approximate solution with the required residual norm.
\end{itemize}

Let $T_{\mathrm{\textit{dp}}}$ denote the average time to perform one CG iteration in double precision, and $T_{\mathrm{\textit{sp}}}$---in single precision; $\nu\!=\!T_{\mathrm{\textit{dp}}}/T_{\mathrm{\textit{sp}}}$ is the speedup achieved by performing computations in single precision instead of double. For example, a server accessible via Google Colab, depending on its load, provides a speedup of $\nu\!\geq\!2$, reaching values of $5$, $7$, and higher.

Let us define the \textit{efficiency} of solving a system of linear equations with a matrix $A$ using a two-stage computation scheme as the quantity $E(A)$, given by the formula $$E(A)=\frac{T_{\mathrm{all\, in\, \textit{dp}}}-T_{\mathrm{mixed}}}{T_{\mathrm{all\, in\, \textit{dp}}}},$$
where $T_{\mathrm{all\, in\, \textit{dp}}}$ is the time required to obtain a solution with a given accuracy when computations are performed in double precision at all iterations, and $T_{\mathrm{mixed}}$ is the time required for computations according to the two-stage scheme.

It is evident that the optimal moment for switching the precision mode (for transitioning to the second stage of the computation scheme) is the iteration preceding $N_{\mathrm{div}}$, at which the divergence of residuals has not yet occurred. By choosing $N_{\mathrm{\textit{sp}}}$ in this way, the maximum number of iterations in double precision is replaced by iterations in single precision, while the required final computational accuracy is achieved at the second stage.

As shown above, the value of $N_{\mathrm{div}}$ correlates with $\mathcal{D}$ and $\tilde{\mathcal{D}}$. For the types of matrices considered by us, the values of $N_{\mathrm{div}}$ turn out to be close if the values of $\mathcal{D}$ and $\tilde{\mathcal{D}}$ are close. Based on this, the choice of the optimal (in terms of efficiency $E(A)$) iteration number $N_{\mathrm{\textit{sp}}}$, at which the precision mode switching occurs, can be made based on the values of $N_{\mathrm{div}}$ precomputed for matrices with $(n,m,\tilde{\mathcal{D}})$ values close to those calculated for the input matrix.

In \cite{Prolubnikov2}, the choice of the optimal iteration number $N_{\mathrm{\textit{sp}}}$ for switching is made indirectly. A threshold of accuracy (residual norm) is selected; upon reaching this threshold at the first stage, the transition to the second stage occurs. The later stagnation occurs in single-precision computation mode, that is, the larger $N_{\mathrm{div}}$ is, the higher the accuracy that can be achieved at the first stage. Using the pseudo-diameter to determine $N_{\mathrm{\textit{sp}}}$ (via a modification of the $k$-nearest neighbors method applied to a small training set) allows achieving efficiency values $E(A)$ that are only approximately $3\%$ lower than the efficiency achieved with the optimal choice of $N_{\mathrm{\textit{sp}}}$. The total computation time required to determine $N_{\mathrm{\textit{sp}}}$ for the input matrix does not exceed $1\%$ of the CG computation time in double precision.

The efficiency increases with the matrix dimension for the classes of matrices considered in this work and in \cite{Prolubnikov2}. Even with a rather pessimistic estimate of $\nu\!=\!3$, the average efficiency $E(A)$ across matrix samples exceeds $20\%$ for random matrices and sparse banded matrices of dimension $n\!=\!1000$, as well as for matrices with extended star graphs of dimension $n\!=\!1001$.\\

\section{Preconditioning}

Consider the preconditioned CG method in the following form:
\begin{table}[h!]
\begin{tabular}{lll}

0. & $k=0$; $d^{0}=r^{0}=M^{-1}(b-Ax^{0})$; & --- \textit{initialization}\\
\quad & \textbf{for} $k\!:=1$ \textbf{to} $N$\\
1. &\qquad $\alpha_{k}=\bigl (r^{k},r^{k}\bigr )/\bigl (d^{k}, Ad^{k}\bigr )$; & --- \textit{step size} \\
2. &\qquad $x^{k+1}=x^{k}+\alpha_{k}d^{k}$; & --- \textit{approximate solution}\\
3. &\qquad $r^{k+1}=r^{k}-\alpha_{k}Ad^{k}$; & --- \textit{recursively computed residual}\\
4. &\qquad $\beta_{k+1}=\bigl (r^{k+1},M^{-1}r^{k+1}\bigr )/\bigl (r^{k},M^{-1}r^{k}\bigr )$; & \\
5. &\qquad $d^{k+1}=M^{-1}r^{k+1}+\beta_{k+1}d^{k}$ & --- \textit{search direction} \\

\end{tabular}
\end{table}
\smallskip

The question of the influence of the metric characteristics of the matrix graph on the convergence of the preconditioned conjugate gradient method is beyond the scope of this work. We only note that when using preconditioning, the relationship between the onset of CG stagnation and the metric characteristics of the matrix graph is preserved, provided that the preconditioning does not disrupt the information propagation structure determined by the original matrix of the system of linear equations. This requirement is met if the term $M^{-1}r^{k+1}$ in the recurrent expression for $d^{k+1}$ does not introduce additional information propagation beyond that provided by the term $Ad^{k}$ from formula~(4), modified for the preconditioned case: 
$$d^{k+1}=(1+\beta_{k+1})M^{-1}d^{k}-\beta_{k}M^{-1}d^{k-1}-\alpha_{k}Ad^{k}.$$ 
This requirement is satisfied for Jacobi preconditioning. As experiments show, incomplete Cholesky factorization (\textit{IC}(0)) for moderately ill-conditioned Laplacian matrices ($10^3\!\le\!\kappa(A)\!\le\!10^4$) also satisfies this requirement.\\

\section*{Conclusions}

A model of propagation and accumulation of errors caused by round-off errors during CG iterations is proposed, based on their interpretation as iterations of traversing the matrix graph of a system of linear equations. In the matrix graph, vertices correspond to the components of the search direction, residual, and approximate solution vectors, which are computed recurrently at~CG iterations.

It is shown that if the number of performed CG iterations is less than the distance between two vertices in the graph, then round-off errors arising in the computation of one component do not affect the errors in the computation of another and do not become correlated. At early iterations, errors propagate locally between neighboring vertices and become locally correlated. As errors become sufficiently global, loss of orthogonality of residuals (conjugacy of search directions) occurs, and stagnation of the iterative process ensues.

Results of computational experiments are presented for matrices with graphs of a given diameter, for Laplacian matrices of differential equations, and for random matrices. These results demonstrate that metric characteristics of the matrix graph, which are efficiently computable for sparse matrices, can be used in mixed-precision algorithms for solving systems of linear equations, namely---for switching precision modes as round-off errors propagate and accumulate.\\

%

\end{document}